\documentclass[12pt,a4paper]{article}
\usepackage[latin1]{inputenc}
\usepackage[T1]{fontenc}
\usepackage{amsmath,amssymb}
\usepackage{geometry}
\newtheorem{theorem}{Theorem}
\newtheorem{definition}{Definition}
\newtheorem{proposition}{Proposition}
\newtheorem{corollary}{Corollary}
\newtheorem{lemma}{Lemma}
\newtheorem{remark}{Remark}
\newtheorem{example}{Example}
\newcommand\Var{\mathrm {Var}}
\newcommand\Bias{\mathrm{Bias}}
\newcommand\Cov{\mathrm{Cov}}
\newcommand\Card{\mathrm{Card}}
\newcommand\argmin{\mathrm{argmin}}
\newcommand\MSE{\mathrm{MSE}}
\newcommand\IMSE{\mathrm{IMSE}}
\renewcommand{\textbf}[1]{\begingroup\bfseries\mathversion{bold}#1\endgroup}
\usepackage{xcolor}
\usepackage{tikz,tkz-tab}
\usepackage{dsfont}
\usepackage{mathrsfs}
\usepackage{hyperref}
\usepackage{multirow}
\usepackage{float}
\makeatletter
\date{{}}
\def\hlinewd#1{%
\noalign{\ifnum0=`}\fi\hrule \@height #1 %
\futurelet\reserved@a\@xhline}
\makeatother
\title{Trapezoidal rule and sampling designs for the nonparametric estimation of the regression function in models with correlated errors}
\begin{document}
\maketitle
\begin{center}

D. BENELMADANI, {\let\thefootnote\relax\footnote{{\textbf{Contact:} djihad.benelmadani@univ-grenoble-alpes.fr, karim.benhenni@univ-grenoble-alpes.fr, sana.louhichi@univ-grenoble-alpes.fr.}}} K. BENHENNI  and S. LOUHICHI\\[0.2cm]

{\it
Laboratoire Jean Kuntzmann (CNRS 5224), Université Grenoble Alpes, France.
}
\end{center}

\textbf{Abstract:} The problem of estimating the regression function in a fixed design models with correlated observations is considered. 
Such observations are obtained from several experimental units, each of them forms a time series. Based on the trapezoidal rule, we propose a simple kernel estimator  and we  derive the asymptotic expression of its integrated mean squared error IMSE and its asymptotic normality.
The problems of the optimal bandwidth and the optimal design with respect to the asymptotic IMSE are also investigated.
Finally, a simulation study is conducted to study the performance of the new estimator and to compare it with the classical estimator of Gasser and M{\"u}ller  in a finite sample set. In addition, we study the robustness of the optimal design with respect to the misspecification of the autocovariance function.
\bigskip

\section{Introduction} 
A classical problem in Statistics is the nonparametric estimation of the regression function of a response variable $Y$ given an explanatory variable $X$, i.e, estimating the function $g$ defined by $g(t)=\mathbb E(Y|X=t)$, based on the observations of $(X_i,Y_i)_{1 \leq i \leq n}$ which are copies of $(X,Y)$. These observations are often modeled as follows: $Y_i =g(t_i)+\varepsilon_i$ where $g$ is the unknown regression function to be estimated, the $\{t_i, i = 1,\cdots,n\}$ is the sampling design and $\{\varepsilon_i, i = 1,\cdots,n\}$ are centered errors. Typically when $(\varepsilon_i)_i$ are i.i.d. the estimation of $g$ has been extensively investigated by several authors. We mention, among others, the work of Priestly and Chao \cite{Priestly & Chao}, Benedetti \cite{Benedetti} and Gasser and M{\"u}ller \cite{Gasser and Muller 1979, Gasser and Muller 1984}.
However, considering that the observations are independent is not always a realistic assumption. In pharmacokinetics for instance, one wishes to estimate the concentration-time of some injected medicine in the organism, based on the observation of blood tests over a period of time. It is clear that the observations provided from the same individual are correlated. For this reason, we shall investigate in this paper the nonparametric regression estimation problem where the observations are correlated. 

\noindent We consider the so-called fixed design regression model with repeated measurements, i.e., \begin{equation}
Y_j(t_i)=g(t_i)+\varepsilon_j(t_i)~~\text{for}~i=1,\cdots,n~\text{and}~j=1,\cdots,m,\label{modelinintro}
\end{equation}
where $\{\varepsilon_j,j=1,\cdots, m\}$ is a sequence of i.i.d. centered error processes with the same distribution as a process $\varepsilon$. Such models are well known in growth curve analysis and in dose response curves. They can be obtained, as noted by Azzalini \cite{Azzalini}, from $m$ individual being observed on a period of time. Generally, observations between different individuals will be uncorrelated. Hence, it is of interest to relax the assumption of correlation between the experimental units.

M{\"u}ller  \cite{Muller 1984} considered Model \eqref{modelinintro} for $m=1$ (observations on one experimental unit) and he supposed that, for $s \ne t$, the covariance $\Cov (\varepsilon_j(t),\varepsilon_j(s))$ tends to $0$ as $n$ tends to infinity, which is not a realistic assumption, as indicated by Hart and Wherly \cite{Hart and Wherly}, in the growth curve problems. They investigated the estimation of $g$ in Model \eqref{modelinintro} with  a stationary error process. They used the estimator proposed by Gasser and M{\"u}ller \cite{Gasser and Muller 1979}, and they showed that, in order to obtain the consistency of the kernel estimator in the presence of correlations, it is necessary to take $m$ experimental units and to let $m$ tends to infinity. 

The stationarity assumption is however restrictive, for instance, in the previous pharmacokinetics example, it is clear that the concentration of the medicine will be high at the beginning then decreases with time. For this, we shall investigate the estimation of $g$ in Model \eqref{modelinintro} where $\varepsilon$ is a nonstationary error process. This case was partially investigated by Ferreira et al.  \cite{Ferreira} and  Benhenni and Rachdi \cite{Benhenni and Rachdi}, where the Gasser and M{\"u}ller estimator was used. 

In this paper, we propose a new estimator for the regression function $g$ as an approximation of the kernel estimator based on continuous observations in the whole interval $[0,1]$ constructed through a stochastic integral. See, for instance, Blanke and Bosq \cite{Blanke and Bosq}, Didi and Louani \cite{Didi and Louani}. When only discrete observations are available, we use the "best" approximation of the stochastic integral, which is obtained by using the trapezoidal rule  based on discrete observations at appropriate $n$ sampling points generated by a sampling density in the interval $[0,1]$.

This estimator has a relatively simpler expression than the kernel estimator proposed by  Gasser and M{\"u}ller \cite{Gasser and Muller 1979}.  Moreover, since this last one depends on $n$ integrals of a kernel at middle samples; and may be subject to numerical (computational) instability,  for instance when a Gaussian kernel is used, whereas the proposed estimator depends only on the observations and the values of the kernel at the sampling points.     
   
In addition to its simple expression, the proposed estimator allows to bring an answer to another important and open statistical problem under correlated errors, which is the optimal design problem. For instance, in the previous pharmacokinetic example, one wishes to find the best moments for the blood testing to be made in order to have a better estimate of the concentration curve.

The optimal design problem has been extensively studied in parametric regression. We mention the work of Sacks and Ylvisaker \cite{Sacks and Ylvisaker}, Belouni and Benhenni \cite{Belouni}  and more recently Dette \textit{et al.} \cite{Dette et al.} among others. In the nonparametric case, M{\"u}ller \cite{Muller 1984} introduced the optimal design points when the errors are asymptotically independent. He used a regular design sequence generated by a density function $f$, i.e, $t_i= F^{-1}(\frac{i}{n})$, where $F$ is the distribution function associated to $f$. He derived the optimal  design generated by a density that minimizes  the asymptotic Integrated Mean Squared Error (IMSE). To the best of our knowledge, there exists no result concerning  the problem of optimal design for nonparametric regression estimation in models under more general class of error processes.

We also investigate the problem of the asymptotic optimal bandwidth. We mention, for the nonparametric case, the work of Hart and Wherly \cite{Hart and Wherly} and Benhenni and Rachdi \cite{Benhenni and Rachdi}. For results on the break down of some data based methods for bandwidth selection in the presence of correlation, for instance the cross validation, and other alternative methods, the reader is referred to Chiu \cite{Chiu 1989}, Altman \cite{Altman 1990},  Hart \cite{Hart 1991, Hart 1994} among others. 

This article is organized as follows. In Section 2, we present the new estimator of the regression function $g$ in Model \eqref{modelinintro} where  $\varepsilon$ is a  centered error process.  In Section 3, we give the asymptotic expressions of the bias, the variance and the IMSE. We then derive the asymptotic optimal bandwidth with respect to the asymptotic IMSE. In addition, we obtain the optimal design density with respect to the asymptotic IMSE, and we prove that it is minimax optimal. We also prove the asymptotic normality of the proposed estimator. In Section 4, we conduct a simulation study to investigate the performance of the new estimator and then to compare it with that of Gasser and M{\"u}ller \cite{Gasser and Muller 1979}. We also conducted a study to compare the uniform and the optimal sampling designs, and to study the robustness of the optimal design, with respect to the misspecification of the autocovariance function.
Since the classical cross validation criteria turned out to be inefficient  in the presence of correlation, we use the bandwidth that minimizes the exact IMSE, the comparison is performed for different numbers of experimental units and  design points.
Finally, Section 5 is dedicated to the proofs of our theoretical results. 

\section{Model and estimator}
We consider $m$ experimental units, each of them having $n$ different measurements of the response (say $ 0 \leq t_1 < t_2 < \dots < t_n \leq 1$). The so-called fixed design regression model is defined as follows:
\begin{equation}
Y_j(t_i)=g(t_i)+\varepsilon_j(t_i) \text{ where } j = 1,\dots, m \text{  and  } i = 1,\dots,n,
\label{model}
\end{equation}
where $g$ is the unknown regression function on $[0,1]$ and $\{\varepsilon_j(t), t \in [0,1] \}_j$ is a sequence of error processes.

We assume that $g \in C^2([0,1])$ and that $(\varepsilon_j )_{j}$ are i.i.d. processes with the same distribution  as a centered second order process $\varepsilon$. We denote by $R$ its autocovariance function.
\subsection{Simple estimator and sampling design}
In order to motivate the construction of our new estimator, we consider the regression model using $m$ continuous experimental units, i.e, 
 \begin{equation}\label{modelcontinue}
Y_j(t)=g(t)+\varepsilon_j(t)~\text{for}~t \in [0,1]~~\text{and } j=1,\cdots,m.
\end{equation}
A continuous kernel estimator of $g$ in Model $\eqref{modelcontinue}$ is given for any $x \in [0,1]$ by, \begin{equation}
\hat{g}_{[0,1]}(x)=\int_0^1 \varphi_{x,h}(t)\overline{Y}(t)~dt~~~\text{with}~~\overline{Y}(t)=\frac{1}{m}\sum_{j=1}^{m}Y_j(t), \label{continuousestimator}
\end{equation}
where $\varphi_{x,h}(t)=\frac{1}{h}K\Big( \frac{x-t}{h}\Big)$ for a kernel $K$ and a bandwidth $h$.
For details on the Kernel estimation of the regression function based on continuous observations see, for instance, Blanke and Bosq \cite{Blanke and Bosq} or  Didi and Louani \cite{Didi and Louani}.

In the practical case where we only have access to discrete observations,  we apply the trapezoidal rule to approximate the continuous Kernel estimator given by \eqref{continuousestimator}. We construct then a new simple estimator of the regression function that we shall call the trapezoidal estimator.

Before introducing the proposed estimator, we begin with defining a sequence of designs which will be used in its construction. This class of designs was considered by Sacks and Ylvisaker \cite{Sacks and Ylvisaker III}.
\begin{definition} \label{regdesigndef}
Let $F$ be a distribution function of some density $f$ satisfying $\underset{t \in [0,1]}{\inf} f(t) > 0$ and $\underset{t \in [0,1]}{\sup} f(t) < \infty$. The so-called regular sequence of designs generated by a density $f$ is defined by,
 \begin{equation*}
 T_n = \bigg \{t_{i,n}=F^{-1}\bigg(\frac{i}{n}\bigg), ~i=1,\dots,n. \bigg\}~~\text{for}~n\geq 1.
 \end{equation*}
 \end{definition} 
 
 Such a sequence of designs verifies the next useful lemma. 
 \begin{lemma}\label{N_T=Onh}
For $n\geq 1$ let $T_n=(t_{i,n})_{i=1,\cdots,n}$ be a regular sequence of designs generated by some density function. Let $x \in ]0,1[$, $h >0$ and note by $N_{T_n} \overset{\Delta}{=} \Card~( T_n \cap [x-h,x+h] )$. Suppose that $N_{T_n}  \ne 0$ and that $nh \geq 1$. Then, \begin{equation}
\underset{0\leq j \leq n}{\sup}\;(t_{j+1,n}- t_{j,n}) =O\Big(\frac{1}{n}\Big)~~\text{ and }~~N_{T_n}=O(nh).\label{cardinal=O(nh)}
\end{equation}
\end{lemma}

We shall now give the definition of the trapezoidal estimator, obtained from  a discrete approximation of the continuous estimator $\hat{g}_{[0,1]}$ given by \eqref{continuousestimator}.
\begin{definition}
The trapezoidal estimator of the regression function $g$ based on the observations $(t_{i,n},Y_j(t_{i,n}))_{\underset{1 \leq j \leq m}{1 \leq i \leq n}}$, where $T_n=(t_{i,n})_{1 \leq i \leq n}$ is a regular sequence of designs generated by a  density function $f$ of support intersecting $[x-h,x+h]$ is given, for any $x \in [0,1]$, by,\begin{equation}\label{trapezestimator}
\hat{g}^{\text{trap}}_n(x) = \frac{1}{2n} \sum_{k=1}^{N_{T_n}-1}\bigg\{ \bigg( \frac{\varphi_{x,h}}{f}\overline{Y}\bigg) (t_{x,k})+\bigg( \frac{\varphi_{x,h}}{f}\overline{Y}\bigg) (t_{x,k+1})\bigg\},\end{equation}
where  $t_{x,1} < \cdots < t_{x,N_{T_n}}$ are the points of $T_n$ in $ [x-h,x+h]$,  $\varphi_{x,h}(t) = \frac{1}{h}K\big(\frac{x-t}{h} \big)$,  $\overline{Y}$ is given in \eqref{continuousestimator},  $K$ is a kernel of support $[-1,1]$ and $h=h(n,m)$ is a bandwidth with $0<h<1$. 
\end{definition}

In order to derive our asymptotic results, the following assumptions on the autocovariance function $R$ and the kernel $K$ are required.  
\subsection{Assumptions}
\begin{itemize}
\item [(A)] The autocovariance function $R$ exists and is continuous 
on the square $[0,1]^2.$
\item [(B)]  At the diagonal (when $t = s$  in the unit square), $R$ has continuous left and right first-order derivatives, that is:   \[
R^{(0,1)}(t,t^-)=\lim\limits_{s\uparrow t}^{}\frac{\partial R(t,s)}{   \partial s} \quad  \text{and} \quad  R^{(0,1)}(t,t^+)=\lim\limits_{s\downarrow t}\frac{\partial R(t,s)}{   \partial s}.
\]
 The jump function along the diagonal $\alpha(t) \overset{\Delta}{=} R^{(0,1)}(t,t^-)-R^{(0,1)}(t,t^+)$ is assumed to be continuous and not identically equal to zero.
\item [(C)]   Off the diagonal (when $t \ne s$  in the unit square), $R$ is assumed to have continuous mixed partial derivatives up to order two and, 
\[A^{(i,j)}\overset{\Delta}{=}\underset{0 \leq t\ne s \leq 1}{\text{ sup }}|R^{(i,j)}(t,s)| < \infty \text{ for }  i,j \text{ such that } 0 \leq i+j \leq 2. \]
\item [(D)] The Kernel $K$ is even at least in $ C^2([-1,1])$ and $K''$ is Lipschitz on [-1,1].
\end{itemize}
Examples of processes with autocovariances satisfying Assumptions $(A),(B)$ and $(C)$ are given as follows.
\begin{example}~~~~~ \label{example1}
\begin{enumerate}
\item The Wiener process with autocovariance function {\small$R(s,t)=\sigma^2\min (s,t)$}, has a constant jump function $\alpha(t)=\sigma^2$ and $R^{(i,j)}(s,t)=0$ for all $i,j$ such that $i+j=2$ and $s\ne t$.
\item The Ornstein-Uhlenbeck process with a stationary autcovariance  $R(s,t)=\sigma^2\exp(-\lambda|s-t|)$ for $\sigma > 0$ and $\lambda > 0$. For this process $\alpha(t)=2\sigma^2\lambda$ and $R^{(0,2)}(s,t)=\sigma^2\lambda^2\exp(-\lambda|s-t|)$.
\item A generalization of the Ornstein-Uhlenbeck process to a process with a nonstationary autocovariance function  of the form: $ R(s,t)=\sigma^2\rho^{|s^\lambda-t^\lambda|/\lambda}$ for $\sigma > 0$, $\lambda > 0$ and $0 < \rho < 1$. For this process the jump function, which is not constant when $\lambda \ne 1$, is given by $\alpha(t)=-2\sigma^2ln(\rho)t^{\lambda-1}$.
\item Sacks and Ylvisaker \cite {Sacks and Ylvisaker} gave another general class of convex stationary autcovariance functions of the form, 
\[R(s,t) = \int_{0}^{1/|t-s|} (1-\mu |t-s|)p(\mu)\;d\mu,\]
where $p$ is a probability density and $p'$ its derivative are such that, 
\[ \underset{\mu \to \infty}{\lim}{\mu^3 p(\mu)} < \infty, ~~~\text{and}~~ \int_{a}^{\infty} (\mu p'(\mu)+3p(\mu))^2)\mu^6 d\mu < \infty,\]
for some finite constant $a$. For this autocovariance function,  $\alpha(t)=2\int_{0}^{\infty} \mu p(\mu)~d\mu$ for all $t$.
\end{enumerate}
\end{example}
The following kernels satisfy Assumption $(D)$.
\begin{example}~~~~
\begin{enumerate}
\item The Quadratic  kernel defined by $K(u)  = \frac{15}{16} (1-u^2)^2~\mathds{1}_{\{|u|\leq 1 \}}$.
\item The Triweight  kernel defined by $K(u)  = \frac{35}{32} (1-u^2)^3~\mathds{1}_{\{|u|\leq 1 \}}$.
\end{enumerate}
\end{example}
\section{Asymptotic results}
The following propositions give the asymptotic expressions of the bias and the variance of the trapezoidal estimator as defined by \eqref{trapezestimator}.
\begin{proposition}\label{biaslemma}
Suppose that Assumption $(D)$ is satisfied. Moreover assume that  $f \in C^2([0,1])$ and $f'',g''$ are Lipschitz functions on $[0,1]$. If $\underset{n \to \infty}{\lim} h = 0$  and $\underset{n \to \infty}{\lim} nh =\infty$ then for any $x \in ]0,1[$,
\[\Bias(\hat{g}^{trap}_n (x) )= \frac{1}{2} h^2 g''(x)B + o (h^2) + O\Big(\frac{1}{n^3h^3}\Big),\]
where $B = \int_{-1}^{1} t^2 K(t) ~dt $.
\end{proposition}
\begin{proposition}
\label{variancelemma1}
Suppose that  Assumptions $(A),(B),(C)$ and $(D)$ are satisfied. Moreover assume that $f \in C^2([0,1])$ and for any $t \in [0,1]$, $f''$ and $R^{(0,2)}(t,.)$  are all Lipschitz on $[0,1]$. If  $\underset{n \to \infty}{\lim} h = 0$  and $\underset{n \to \infty}{\lim} nh =\infty$ then for any $x \in ]0,1[$,
\begin{align*}
\Var (\hat{g}^{trap}_n (x) ) & = \frac{1}{m} \Big (R(x,x)-\frac{h}{2}C_K \alpha (x) \Big )+\frac{V}{12mn^2h}\frac{\alpha(x) }{f^2(x)}\\
&~~~~~ + o\Big(\frac{h}{m}\Big)+O\Big(\frac{1}{mn^2}+\frac{1}{mn^3h^3}\Big),
\end{align*}
where $V = \int_{-1}^{1} K^2(t)\;dt$ and $C_K = \int_{-1}^1 \int_{-1}^1 |u-v|K(u)K(v) du dv.$
\end{proposition}

Propositions \ref{biaslemma} and \ref{variancelemma1} allow to derive the asymptotic expression of the mean squared error (MSE) of the Trapezoidal  estimator \eqref{trapezestimator}. The  integrated mean squared error (IMSE) is then obtained by integrating the MSE with respect to some weight function $w$. The results are announced, without proof, in the following theorem.
\begin{theorem}\label{IMSE}
If all the assumptions of Propositions  \ref{biaslemma} and \ref{variancelemma1} are satisfied then for any $x \in ]0,1[$,
{\small
\begin{align*}
\MSE({\hat{g}^{trap}_n} (x)) & = \frac{1}{m}\Big ( R(x,x) -\frac{h}{2}\alpha(x) C_K \Big)+\frac{V}{12mn^2h}\frac{\alpha(x) }{f^2(x)}+ \frac{1}{4} h^4 [g''(x)]^2 B^2\\
&~~~~ +o \big(h^4+\frac{h}{m}\big)+O\Big(\frac{1}{n^3h}+\frac{1}{mn^3h^3}+\frac{1}{mn^2}+\frac{1}{n^6h^6}\Big),
\end{align*}
}
{\small
\begin{align}
\IMSE ({\hat{g}^{trap}_n}) & = \frac{1}{m} \int_{0}^{1}\Big ( R(x,x) -\frac{h}{2}\alpha(x) C_K \Big) w(x)\;dx + \frac{V}{12mn^2h} \int_{0}^{1}  \frac{\alpha(x)}{f^2(x)} w(x)\;dx \nonumber \\
&~~~~~~ + \frac{1}{4} h^4 B^2 \int_{0}^{1}[g''(x)]^2 w(x)\;dx +o \big(h^4+\frac{h}{m}\big)\nonumber \\
&~~~~~~+O\Big(\frac{1}{n^3h}+\frac{1}{mn^3h^3}+\frac{1}{mn^2}+\frac{1}{n^6h^6}\Big), \label{IMSEexpression}
\end{align} }
where $w$ is a continuous density function, $V$, $B$ and $C_K$ are given in Propositions \ref{biaslemma}, \ref{variancelemma1}.
\end{theorem}
The previous Theorem shows, the efficiency of the Trapezoidal estimator, since the IMSE tends to $0$ when $m \to \infty$, $h \to 0$ and $nh \to \infty$ as $n \to \infty$.

The asymptotic optimal bandwidth is obtained by minimizing the asymptotic IMSE as given by the following proposition.
\begin{proposition}[Optimal bandwidth]\label{optimalh}
Suppose that the assumptions of Theorem \ref{IMSE} are satisfied. Moreover assume that $  \frac{m}{n} = O(1)$ as $n,m \to \infty$. Denote by IMSE($h$) the IMSE of the trapezoidal estimator when the bandwidth $h$ is used. Then the bandwidth, 
\begin{equation} \label{hoptimal}
 h^* = \Bigg( \frac{C_K\int_0^1 \alpha(x)w(x)\;dx}{2B^2\int_0^1 [g''(x)]^2 w(x)\;dx}\Bigg)^{1/3}m^{-1/3},
 \end{equation}  is optimal in the sense that, \[\underset{n,m \to \infty}{\overline{\lim}}~\frac{\IMSE(h^*)}{\IMSE(h_{n,m})} \leq 1,\]
for any sequence of bandwidths $h_{n,m}$ verifying: \[\underset{n,m \to \infty}{\lim} h_{n,m} = 0 ~~~\text{and}~~\underset{n,m \to \infty}{\overline{\lim}}~mh_{n,m}^3 < +\infty,\]
where $B$ and $C_K$ are given in Propositions  \ref{biaslemma} and \ref{variancelemma1}.
 \end{proposition}
We are interested now in finding the optimal design density, i.e, a function $f^*$ according to the criteria $f^*\in\underset{f}{\argmin}$ IMSE, where the minimum is taken with respect to the class of positive densities defined on $[0,1]$. In view of Theorem \ref{IMSE}, the asymptotic optimal design density verifies, \[f^*\in\underset{f>0,~\int_0^1 f(x)dx=1}{\argmin} \int_{0}^{1}\frac{\alpha(x) }{f^2(x)} w(x)~dx.\]
This optimization problem is solved in the following corollary. 
\begin{corollary}[Optimal design]\label{optimaldensitytheorem}
Suppose that the assumptions of Theorem \ref{IMSE} are satisfied. If $\underset{n \to \infty}{\lim} nh^2 = \infty$ and $ \underset{n,m \to \infty}{\lim} \frac{n}{m} = \infty$, then the optimal sampling density with respect to the asymptotic $\IMSE$ is given by,
\begin{equation}
f^* (t) = \frac{\{\alpha(t)w(t)\}^{1/3}}{\int_{0}^{1}\{\alpha(s)w(s)\}^{1/3}\;ds}1_{[0,1]}(t). \label{f^*}
\end{equation}
\end{corollary}
Let $\hat{g}^{trap}_{n,f^*}$ be the Trapezoidal estimator \eqref{trapezestimator} with $f=f^*$ defined by  \eqref{f^*}. We have,
\begin{align*}
\IMSE ({\hat{g}^{trap}_{n,f^*}}) & = \frac{1}{m} \int_{0}^{1}\Big ( R(x,x) -\frac{1}{2}\alpha(x) C_K h\Big) w(x)\;dx\\
&~~~~~~+ \frac{V}{12mn^2h} \bigg(\int_{0}^{1} (\alpha(x) w(x))^{1/3}\;dx\bigg)^3+ \frac{1}{4} h^4 B^2 \int_{0}^{1}[g''(x)]^2 w(x)\;dx\\
&~~~~~~+o \big(h^4+\frac{h}{m}\big)+O\Big(\frac{1}{n^3h}+\frac{1}{mn^3h^3}+\frac{1}{mn^2}+\frac{1}{n^6h^6}\Big).
\end{align*}
\begin{remark}
Let  $\hat{g}^{trap}_{n,unif}$ be the Trapezoidal estimator \eqref{trapezestimator} with a uniform density, i.e, $f=f_{unif}$ the identity in $[0,1]$. The asymptotic IMSE of  $\hat{g}^{trap}_{n,unif}$ is given by,
\begin{align*}
\IMSE ({\hat{g}^{trap}_{n,unif}}) & = \frac{1}{m} \int_{0}^{1}\Big ( R(x,x) -\frac{1}{2}\alpha(x) C_K h\Big) w(x)\;dx + \frac{V}{12mn^2h} \int_{0}^{1}  \alpha(x) w(x)\;dx \nonumber \\
&~~~~~~+ \frac{1}{4} h^4 B^2 \int_{0}^{1}[g''(x)]^2 w(x)\;dx+o \big(h^4+\frac{h}{m}\big)\nonumber \\
&~~~~~~+O\Big(\frac{1}{n^3h}+\frac{1}{mn^3h^3}+\frac{1}{mn^2}+\frac{1}{n^6h^6}\Big).
\end{align*} 
The reduction of the residual IMSE, $\overline{\IMSE} \overset{\Delta}{=} \IMSE -\sigma_{x,h}^2/m$, by using the asymptotic optimal design over the uniform design is then,
\begin{align*}
rIMSE=\frac{\overline{\IMSE} ({\hat{g}^{trap}_{n,unif}})-\overline{\IMSE} ({\hat{g}^{trap}_{n,f^*}})}{\overline{\IMSE} ({\hat{g}^{trap}_{n,unif}})} & \sim 1 - \frac{\big(\int_{0}^{1}(\alpha(x) w(x))^{1/3}\;dx\big)^3}{\int_{0}^{1}  \alpha(x) w(x)\;dx}.
\end{align*}
For instance, if $R(s,t)=st \min(s,t)$ then $\alpha(t)=t^2$. Taking $w \equiv 1$ gives $rIMSE \sim  35\%$.
\end{remark}
Finally, the next theorem gives the asymptotic normality of the Trapezoidal estimator \eqref{trapezestimator}.

The optimal design, generated by the density function \eqref{f^*}, may not be robust with respect to the missspecification of the autocovariance jump function $\alpha$, and the weight function $w$. For this, we shall use a minimax criterion to obtain the optimal sampling design. Biedermann and Dette \cite{Biedremann} gave the following criterion, a density function $f^*$ is said to be minimax optimal if,
\begin{equation}\label{minimaxcrit}
f^*\in \underset{f>0, \int_{0}^{1}f(t)dt=1}{\argmin}~\underset{(\alpha, w) \in \Lambda}{\max} \Psi_{(\alpha, w)} (f),
\end{equation}
where, $$\Psi_{(\alpha, w)} (f) = \int_{0}^{1} \frac{\alpha(t)}{f^2(t)}w(t)dt,$$
and, $$ \Lambda = \Big\{(\alpha, w) \in (C[0,1])^2 ~\big / \int_0^1 \alpha(t)dt < \epsilon_1, \bigg( \int_0^1 w(s)^{1/2}ds\bigg)^2< \epsilon_2 \Big\}.$$
The following theorem assures that the asymptotic optimal design density, defined in Corollary \ref{optimaldensitytheorem}, is optimal in the sense of minimax. 
\begin{theorem}[Minimax optimality]\label{minimaxtheorem}
Suppose that the assumptions of Theorem \ref{IMSE} are satisfied.   The function $f^*$ given by \eqref{f^*} is optimal with respect to the minimax criterion \eqref{minimaxcrit}.
\end{theorem}
Finally, we conclude our theoretical results by the asymptotic normality of the trapezoidal estimator, presented in the following theorem.
\begin{theorem}[Asymptotic normality]\label{asymptotic normality theorem}
Suppose that the assumptions of Theorem \ref{IMSE} are satisfied. If $\underset{m \to \infty}{\lim}\sqrt mh^2 = 0$ and $\underset{n \to \infty}{\lim} nh^2 = \infty$ then for any $x \in ]0,1[$,
\[ \sqrt m \Big( {\hat{g}^{trap}_n}(x) - g(x) \Big) \overset{\mathscr{D}}{\longrightarrow} Z,~~~\text{with}~Z \sim \mathcal N (0,R(x,x)), \]
where $\mathscr D$ denotes the convergence in distribution and $\mathcal N$ is the normal distribution.
\end{theorem}
\section{Simulation study}
In this section, we  investigate the performance of our estimator \eqref{trapezestimator} in a finite sample set.  We shall use the  cubic growth curve, used by Benhenni and Rachdi \cite{Benhenni and Rachdi} and Hart and Wherly \cite{Hart and Wherly}, 
\begin{equation}
\label{cubic}
g(x)=10x^3-15x^4+6x^5~~~\text{for}~~0<x<1.
\end{equation}
This function was mainly used  due to its similarity to the logistic function which is frequently found in growth curve analysis.
The sampling points are taken to be:  
\begin{equation}
t_i=(i-0.5)/n~~\text{for}~i=1,\cdots,n. \label{designinsimulation}
\end{equation}
The error process $\varepsilon$ is taken to be the Wiener error process with autocovariance function $R(s,t)=\sigma^2\min(s,t)$. 
The Kernel used here is the quadratic kernel given by $K(u)=(15/16)(1-u^2)^2I_{[-1,1]}(u)$. The bandwidth used in this study is the optimal bandwidth with respect to the exact $\IMSE$.

\begin{figure}[H] 
\begin{minipage}[c]{.3\linewidth}
     \begin{center}
             \includegraphics[scale=0.49]{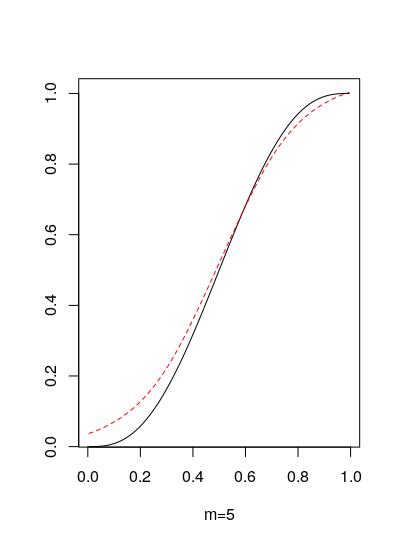}  
         \end{center}
   \end{minipage} \hfill
   \begin{minipage}[c]{.3\linewidth}
     \begin{center}
             \includegraphics[scale=0.49]{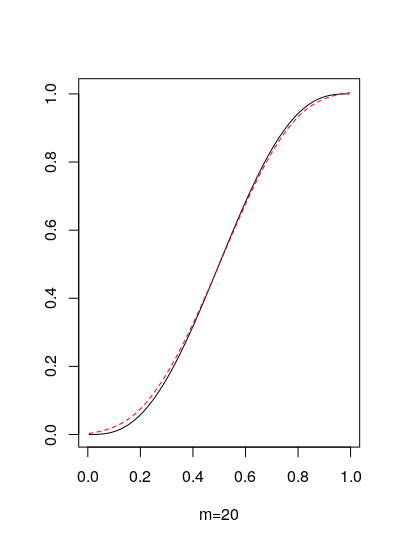}
         \end{center}
   \end{minipage} \hfill
   \begin{minipage}[c]{.3\linewidth}
    \begin{center}
            \includegraphics[scale=0.49]{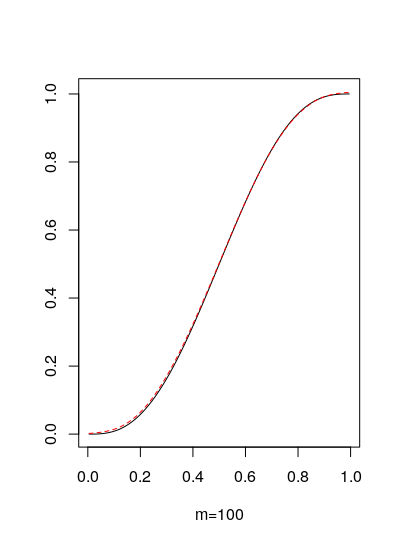}
        \end{center}
 \end{minipage}
\caption{Cubic regression function is in plain line and the trapezoidal estimator is in dashed one.} 
\label{graphsdetrapeze}
\end{figure}

 We consider the mean of all estimators obtained from 100 simulations. We take $\sigma^2=0.5$ and simulations for other values of $\sigma^2$ gave similar results, they are given in Figure \ref{graphsdetrapeze}  for  a fixed number of observations $n=100$ and three different values of experimental units $m=5,20,100$. 
 It is clear that, the performance of the trapezoidal estimator  gets better as $m$ increases.
 
 Our aim now is to compare the trapezoidal estimator to that of Gasser and M{\"u}ller  \cite{Gasser and Muller 1979} (referred by GM estimator), given  for any $x \in ]0,1[$ by,
\begin{equation}\label{GMestimator}
\hat{g}^{GM}_n(x) = \sum_{i=1}^n \int_{m_{i-1}}^{m_i} \varphi_{x,h}(t)dt ~\overline{Y}(t_i),
\end{equation}
where $m_{0}=0$, $m_n=1$ and $m_{i}=(t_i+t_{i+1})/2$ for $i=2,\cdots,n-1$,  $\varphi_{x,h}(t)=(1/h)K((x-t)/h)$ and $\overline{Y}(t_i)=(1/m)\sum_{j=1}^m Y_j(t_i)$.

This comparison is conducted with respect to the non-asymptotic IMSE and under different types of correlation errors. We consider again the cubic regression function, the design given by \eqref{designinsimulation} and the quadratic kernel.  The two error processes considered here are the stationary Ornstein-Uhlenbeck process with $R(s,t)= \exp(-\lambda |s-t|)$, and the nonstationary  Wiener process  with $R(s,t)=\sigma^2 \min(s,t)$. We investigate various "amount" of correlation by taking different values of both $\sigma^2$ and $\lambda$.

We take the weight density $w$ to be uniform on $[0,1]$, and we  compare the optimal non-asymptotic $\IMSE$ of the two estimators, i.e., $\underset{0 < h< 1}{\inf}~ \IMSE (h)$. The bandwidth $h$ is chosen over a grid from $0.09$ to $0.5$. The results are given in Tables \ref{Wiener sigma=1}-\ref{OU sigma=1 alpha=50} for  $n=30$  and for different values of $m$. The tables present the integrated bias squared denoted by $Ibias^2$, integrated variance denoted by  $Ivar$ and the $\IMSE$ together with the optimal bandwidth  associated to the smallest non-asymptotic IMSE for each estimator. The tables are organized according to the "degree" of correlation of the errors.

It can be seen that the optimal bandwidth is the same for both estimators, in addition, as expected, it decreases as  $m$ increases.  

Consider first the case of strong correlated errors, i.e, for a large $\sigma^2$ and a small $\lambda$. In 
Table \ref{Wiener sigma=1}, for the Wiener process with $\sigma^2=1$, it appears that the G-M estimator has a slightly smaller $Ibias^2$  while the trapezoidal estimator has a slightly smaller $Ivar$ and since the $Ibias^2$  is too small compared to the  $Ivar$ then the trapezoidal estimator has a slightly smaller IMSE. 
For the Ornstein-Uhlenbeck with  $\lambda=1$ (c.f. Table \ref{OU sigma=1 alpha=1}) it can be seen that the trapezoidal estimator has a slightly better performance because of a smaller IMSE,  due to a smaller $Ibias^2$  and a smaller $Ivar$.

Consider now the case of moderate correlated errors. In  Table \ref{Wiener sigma=0.5} (for the Wiener process with $\sigma^2=0.5$) it seems that the G-M estimator has a slightly smaller $Ibias^2$  while the trapezoidal estimator has a slightly smaller $Ivar$ and smaller IMSE. While for the Ornstein-Uhlenbeck process with $\lambda=25$, presented in table \ref{OU sigma=1 alpha=25}, the G-M estimator has slightly smaller IMSE  due to a smaller $Ibias^2$  and a smaller $Ivar$.

Finally,  consider the weakly correlated errors, i.e, for a small value of $\sigma^2$ and a large value of $\lambda$. In table \ref{Wiener sigma=0.06}, for the Wiener process with $\sigma^2=0.06$. it appears that the G-M estimator has a slightly smaller $Ibias^2$  while the trapezoidal estimator has a smaller $Ivar$ and smaller IMSE.  However, for the Ornstein-Uhlenbeck process with $\lambda=50$ (c.f. Table \ref{OU sigma=1 alpha=50}) the trapezoidal estimator has a slightly smaller $Ibias^2$  while the G-M estimator has a slightly smaller $Ivar$ and IMSE.

Overall, the two estimators, i.e, the trapezoidal estimator and the Gasser and M{\"u}ller estimator, have "approximately" the same performance. Hence, the proposed estimator, which has a simpler expression, is as efficient as the classical Gasser and M{\"u}ller estimator.

In all the previous cases, it appears that $Ibias^2$ is always smaller than $Ivar$.
It should be noted here that, both of the estimators have boundary problems. A modified kernel at the edges, as suggested by Hart and Wherly \cite{Hart and Wherly}, was used in this simulation. 
\subsection{Optimal design}  
Another important aspect we looked at in this simulation study was the use of the asymptotic optimal design in a finite sample set. We consider the class of autocovariance functions introduced in Example \ref{example1} as follows: $$R(s,t)=\sigma^2 \rho^{|s^\lambda-t^\lambda|/\lambda}, ~~\sigma^2>0, \lambda>0~\text{and}~~0<\rho<1,$$
for which the jump function $\alpha(t)=-2\sigma^2 ln(\rho)t^{\lambda-1}$.
In particular  when $\lambda = 1$ we obtain an Ornstein-Uhlenbeck  stationary error process. In our numerical studies we will consider the nonstationary case,  $\lambda \ne 1$. This class of nonstationary parametric autocovariance functions was introduced by  Núñez-Antón and Woodworth  \cite{Núñez} to study the efficacy of cochlear implants. It was also used by several other authors, by Ferreira et al. \cite{Ferreira} who were interested in obtaining the optimal bandwidth for the Gasser and Müller estimator, by Ziemmerman et al. \cite{Zimmerman}, and then by Núñez-Antón \cite{Núñez 1997} to study the speech recognition data.

We compare, for $m \in \{5, 10, 20, 30\}$ and for instance $h=0.123$, the non-asymptotic IMSE (taking $w \equiv 1$) of the trapezoidal  estimator \eqref{trapezestimator},  using both the uniform design \eqref{designinsimulation}, i.e.,  $f \equiv 1$ and the optimal design generated by $f^*$ given in \eqref{f^*}, i.e., \[f_{\lambda}^*(t)=\frac{\lambda + 2}{3}~t^{(\lambda -1)/3}~1_{[0,1]}(t)~~ \text{and} ~~t_{\lambda,i}^{*}=\Big(\frac{i}{n}\Big)^{3/(\lambda+2)}.\]
\subsubsection*{Robustness of the optimal design}
The optimal design depends on the autocovariance parameter $\lambda$, which is not known in practice, therefore we cannot use this design to compute the estimate of  the regression function $g$. 
As an alternative, we can estimate first the autocovariance parameter $\lambda$, from the observations obtained following a uniform design, then we  obtain the estimated optimal design defined as follows:
\[f_{\widehat{\lambda}}^*(t)=\frac{\widehat{\lambda} + 2}{3}~t^{(\widehat{\lambda} -1)/3}~1_{[0,1]}(t)~~ \text{and} ~~t_{\widehat{\lambda},i}^{*}=\Big(\frac{i}{n}\Big)^{3/(\widehat{\lambda}+2)}.\]
The estimation of the autocovariance parameters is obtained by minimizing the following criterion, as done for instance in Ferreira et al. \cite{Ferreira}: 
\begin{align}
Q_{n,m}(\sigma^2,\lambda,\rho)=\frac{1}{n^2} \sum_{i=1}^n \sum_{j=1}^n \Big(\widehat{R}(t_i,t_j)-R(t_i,t_j) \Big)^2, \label{Qn}
\end{align}
where the empirical correlation estimator is given as follows: \[\widehat{R}(t_i,t_j) =\frac{1}{m-1} \sum_{k=1}^{m} \big( Y_k(t_i)-\overline{Y}(t_i)\big)\big( Y_k(t_j)-\overline{Y}(t_j)\big)~~\text{for}~i,j=1,\cdots,n.\]
From Amemiya  \cite{Amemiya}, as noted by Ferreira and al.  \cite{Ferreira}, it is known that the non linear least square estimator $(\widehat{\sigma^2},\widehat{\lambda},\widehat{\rho})$ is consistent.

In our simulation study, we fixed $\lambda=4, \sigma^2=0.5$ and $\rho=0.5$.
To estimate $(\lambda, \sigma^2, \rho)$, we generated $100$ matrices $(Y_j(t_i))_{\underset{1 \leq j \leq m}{1 \leq i \leq n}}$ of observations using the uniform design. For every matrix, we used the Generalized Simulated Annealing (GSA) algorithm to minimize the function \eqref{Qn}, the estimation  $(\widehat{\lambda}, \widehat{\sigma^2}, \widehat{\rho})$ is then the median of the 100 estimated values.
For more details on the use of the software R algorithm function, see Xiang et al.  \cite{Xiang}. This algorithm is essentially known for its ability to handle very complex non-linear objective functions with a very large number of optima.

The results are given in Tables \ref{GMuTRAPo5}-\ref{GMuTRAPo30}, where the reduction in the IMSE by taking the optimal design instead of the uniform design is given by,
\begin{align*}
rIMSE_{\lambda}=\frac{\IMSE ({\hat{g}^{trap}_{n,unif}})-\IMSE ({\hat{g}^{trap}_{n,f_{{\lambda}}^*}})}{\IMSE ({\hat{g}^{trap}_{n,unif}})},
\end{align*}
and the reduction in the IMSE by taking the plug-in estimated optimal design instead of the uniform design is given by,
\begin{align*}
{rIMSE}_{\widehat{\lambda}}=\frac{\IMSE ({\hat{g}^{trap}_{n,unif}})-\IMSE ({\hat{g}^{trap}_{n,f_{\widehat{\lambda}}^*}})}{\IMSE ({\hat{g}^{trap}_{n,unif}})}.
\end{align*}
It can be seen in Tables \ref{GMuTRAPo5}-\ref{GMuTRAPo30} that there exists a reduction of the IMSE of the Trapezoidal  estimator when using the optimal design, even for small values of the sampling size $n$ and the number of experimental units $m$. Likewise, the estimated optimal design obtained by estimating the covariance parameter, still provides a reduction of the IMSE over the uniform design. This reduction is close to the one using the theoretical optimal design, this shows that the optimal design is robust when the covariance parameter has to be estimated.

\begin{table}[H]
\caption{The integrated squared bias, Integrated variance, IMSE and the optimal bandwidth in terms of $m$ under the Wiener error process with $\sigma^2=1, $ for the GM and the trapezoidal estimator.} \label{Wiener sigma=1}
\begin{center}
\begin{tabular}{|c|c|c|c|c|c|}
\hline 
$n=20$ &$m$ & $Ibias^2$ & $Ivar$ & $\IMSE$ & $h_{opt}$ \\
\hlinewd{1.5pt}
$GM$ &\multirow{2}{*}{5}& 2.8832$\times 10^{-3}$ & 8.4967$\times 10^{-2}$ &8.7850$\times 10^{-2}$ &0.411\\
\cline{1-1}\cline{3-6}
$Trap$&& 2.8833$\times 10^{-3}$ &8.4959$\times 10^{-2}$ &8.7843$\times 10^{-2}$ &0.411 \\
\hlinewd{1.5pt}
$GM$&\multirow{2}{*}{15} & 1.04816$\times 10^{-3}$ &2.9293$\times 10^{-2}$ &3.0341$\times 10^{-2}$ &0.322\\
\cline{1-1}\cline{3-6}
$Trap$& & 1.04856$\times 10^{-3}$ &2.9276$\times 10^{-2}$ &3.0325$\times 10^{-2}$ &0.322 \\
\hlinewd{1.5pt}
$GM$&\multirow{2}{*}{30} & 2.7691$\times 10^{-4}$ &1.5169$\times 10^{-2}$& 1.5446$\times 10^{-2}$ &0.233\\
\cline{1-1}\cline{3-6}
$Trap$& & 2.8535$\times 10^{-4}$ &1.5124$\times 10^{-2}$ &1.5409$\times 10^{-2}$ &0.233 \\
\hline 
\end{tabular}
\end{center}
\label{tab:multicol}
\end{table}
\begin{table}[H]
\caption{The integrated squared bias, Integrated variance, IMSE and the optimal bandwidth in terms of $m$ under the Ornstein-Uhlenbeck error process with  $\lambda=1$ for the GM and the trapezoidal estimator.}\label{OU sigma=1 alpha=1} 
\begin{center}
\begin{tabular}{|c|c|c|c|c|c|}
\hline 
$n=20$ &$m$& $Ibias^2$ & $Ivar$ & $\IMSE$ & $h_{opt}$ \\
\hlinewd{1.5pt}
$GM$&\multirow{2}{*}{5}& 4.57002$\times 10^{-3}$ &1.70570$\times 10^{-1}$ &1.75140$\times 10^{-1}$ &0.46\\
\cline{1-1}\cline{3-6}
$Trap$& &  4.57001$\times 10^{-3}$ &1.70565$\times 10^{-1}$& 1.75135$\times 10^{-1}$& 0.46 \\
\hlinewd{1.5pt}
$GM$&\multirow{2}{*}{15} &1.31050$\times 10^{-3}$ &5.8884$\times 10^{-2}$ &6.0194$\times 10^{-2}$ &0.34 \\
\cline{1-1}\cline{3-6}
$Trap$& &1.30997$\times 10^{-3}$ &5.8857$\times 10^{-2}$ &6.0167$\times 10^{-2}$ &0.34 \\
\hlinewd{1.5pt}
$GM$&\multirow{2}{*}{30}&7.7889$\times 10^{-4}$& 2.9818$\times 10^{-2}$ & 3.0597$\times 10^{-2}$ &0.30 \\
\cline{1-1}\cline{3-6}
$Trap$& &7.7828$\times 10^{-4}$ &2.9791$\times 10^{-2}$ & 3.0569$\times 10^{-2}$ &0.30\\
\hline
\end{tabular}
\end{center}
\label{tab:multicol}
\end{table}
\begin{table}[H]
\caption{The integrated squared bias, Integrated variance, IMSE and the optimal bandwidth in terms of $m$ under the  Wiener error process with $\sigma^2=0.5$  for the GM and the trapezoidal estimator.} \label{Wiener sigma=0.5}
\begin{center}
\begin{tabular}{|c|c|c|c|c|c|}
\hline 
$n=20$&$m$ & $Ibias^2$ & $Ivar$ & $\IMSE$ & $h_{opt}$ \\
\hlinewd{1.5pt}
$GM$&\multirow{2}{*}{5} & 1.0481$\times 10^{-3}$ &4.3939$\times 10^{-2}$ &4.4988$\times 10^{-2}$ &0.322 \\
\cline{1-1}\cline{3-6}
$Trap$& &1.0485$\times 10^{-3}$ &4.3915$\times 10^{-2}$ &4.4963$\times 10^{-2}$ &0.322\\
\hlinewd{1.5pt}
$GM$&\multirow{2}{*}{15} &2.7691$\times 10^{-4}$ & 1.5169$\times 10^{-2}$ & 1.5446$\times 10^{-2}$ & 0.233\\
\cline{1-1}\cline{3-6}
$Trap$&& 2.8535$\times 10^{-4}$ & 1.5124$\times 10^{-2}$ & 1.5409$\times 10^{-2}$ & 0.233\\
\hlinewd{1.5pt}
$GM$&\multirow{2}{*}{30} &1.1792$\times 10^{-4}$ &7.7228$\times 10^{-3}$ &7.8407$\times 10^{-3}$ &0.188 \\
\cline{1-1}\cline{3-6}
$Trap$&& 1.4175$\times 10^{-4}$ &7.6733$\times 10^{-3}$ &7.8150$\times 10^{-3}$ &0.188 \\
\hline
\end{tabular}
\end{center}
\label{tab:multicol}
\end{table}

\begin{table}[H]
\caption{The integrated squared bias, Integrated variance, IMSE and the optimal bandwidth in terms of $m$ under the  Ornstein-Uhlenbeck error process with $\lambda=25$ for the GM and the trapezoidal estimator.} \label{OU sigma=1 alpha=25}
\begin{center}
\begin{tabular}{|c|c|c|c|c|c|}
\hline 
$n=20$&$m$ & $Ibias^2$ & $Ivar$ & $\IMSE$ & $h_{opt}$ \\
\hlinewd{1.5pt}
$GM$&\multirow{2}{*}{5} & 4.3931$\times 10^{-3}$ &2.7163$\times 10^{-2}$ &3.1556$\times 10^{-2}$ &0.455\\
\cline{1-1}\cline{3-6}
$Trap$& &4.3930$\times 10^{-3}$ &2.7165$\times 10^{-2}$ &3.1558$\times 10^{-2}$ &0.455\\
\hlinewd{1.5pt}
$GM$&\multirow{2}{*}{15} &1.7942$\times 10^{-3}$ &1.2819$\times 10^{-2}$ &1.4613$\times 10^{-2}$ &0.366\\
\cline{1-1}\cline{3-6}
$Trap$&& 1.7935$\times 10^{-3}$ &1.2824$\times 10^{-2}$ &1.4618$\times 10^{-2}$ &0.366\\
\hlinewd{1.5pt}
$GM$&\multirow{2}{*}{30} &1.0481$\times 10^{-3}$& 7.0808$\times 10^{-3}$& 8.1290$\times 10^{-3}$& 0.322\\
\cline{1-1}\cline{3-6}
$Trap$& & 1.0485$\times 10^{-3}$ & 7.0855$\times 10^{-3}$ & 8.1341$\times 10^{-3}$ &0.322 \\
\hline
\end{tabular}
\end{center}
\label{tab:multicol}
\end{table}

\begin{table}[H]
\caption{The integrated squared bias, Integrated variance, IMSE and the optimal bandwidth in terms of $m$ under the  Wiener error process with $\sigma^2=0.06$ for the GM and the trapezoidal estimator.} \label{Wiener sigma=0.06}
\begin{center}
\begin{tabular}{|c|c|c|c|c|c|}
\hline 
$n=20$&$m$ & $Ibias^2$ & $Ivar$ & $\IMSE$ & $h_{opt}$ \\
\hlinewd{1.5pt}
$GM$&\multirow{2}{*}{5}& 9.9714$\times 10^{-5}$ &5.5781$\times 10^{-3}$& 5.6778$\times 10^{-3}$ & 0.181 \\
\cline{1-1}\cline{3-6}
$Trap$&&1.2841$\times 10^{-4}$& 5.5373$\times 10^{-3}$ &5.6657$\times 10^{-3}$& 0.181\\
\hlinewd{1.5pt}
$GM$&\multirow{2}{*}{15}&9.9714$\times 10^{-5}$& 4.6484$\times 10^{-3}$ & 4.7481$\times 10^{-3}$& 0.181\\
\cline{1-1}\cline{3-6}
$Trap$&&1.2841$\times 10^{-4}$ &4.6145$\times 10^{-3}$ &4.7429$\times 10^{-3}$ & 0.181\\
\hlinewd{1.5pt}
$GM$&\multirow{2}{*}{30} &9.9714$\times 10^{-4}$ &3.9844$\times 10^{-3}$ &4.0841$\times 10^{-3}$ &0.181\\
\cline{1-1}\cline{3-6}
$Trap$& &1.2841$\times 10^{-4}$& 3.9552$\times 10^{-3}$ &4.0836$\times 10^{-3}$ & 0.181\\
\hline
\end{tabular}
\end{center}
\label{tab:multicol}
\end{table}

\begin{table}[H]
\caption{The integrated squared bias, Integrated variance, IMSE and the optimal bandwidth in terms of $m$ under the  Ornstein-Uhlenbeck error process with $\lambda=50$ for the GM and the trapezoidal estimator.} \label{OU sigma=1 alpha=50}
\begin{center}
\begin{tabular}{|c|c|c|c|c|c|}
\hline 
$n=20$ &$m$ & $Ibias^2$ & $Ivar$ & $\IMSE$ & $h_{opt}$ \\
\hlinewd{1.5pt}
$GM$&\multirow{2}{*}{5} &4.3496$\times 10^{-3}$ & 1.9905$\times 10^{-2}$ &2.4255$\times 10^{-2}$ &0.454 \\
\cline{1-1}\cline{3-6}
$Trap$&& 4.3494$\times 10^{-3}$ &1.9907$\times 10^{-2}$ &2.4257$\times 10^{-2}$ &0.454 \\
\hlinewd{1.5pt}
$GM$&\multirow{2}{*}{15}& 2.8194$\times 10^{-3}$ & 1.8049$\times 10^{-2}$ &2.0868$\times 10^{-2}$ &0.408 \\
\cline{1-1}\cline{3-6}
$Trap$&& 2.8192$\times 10^{-3}$ & 1.8053$\times 10^{-2}$& 2.0872$\times 10^{-2}$& 0.408 \\
\hlinewd{1.5pt}
$GM$&\multirow{2}{*}{30} & 2.8194$\times 10^{-3}$ & 1.5470$\times 10^{-2}$ & 1.8290$\times 10^{-2}$ & 0.408 \\
\cline{1-1}\cline{3-6}
$Trap$& &2.8192$\times 10^{-3}$ &1.5474$\times 10^{-2}$ &1.8293$\times 10^{-2}$& 0.408 \\
\hline
\end{tabular}
\end{center}
\label{tab:multicol}
\end{table}
\begin{table}[H]
\caption{The IMSE and the reductions in the IMSE of $\hat{g}_n^{trap}$ using the uniform design, theoretical optimal design and estimated optimal design when $R(s,t)=\sigma^2 \rho^{|s^\lambda-t^\lambda|/\lambda}$ and $n=5.$ }\label{GMuTRAPo5}
\begin{center}
\begin{tabular}{|l|c|r|c|c|c|c|} 
  \hline
  $m$ & $Trap_{unif}$& $Trap_{opt}$& ${rIMSE}_{\lambda}$ & $Trap_{\widehat{opt}}$ & ${rIMSE}_{\widehat{\lambda}}$ & $\widehat{\lambda}$\\
  \hline
  5 & 0.3661 & 0.3138 & 14.28\%& 0.3167 & 13.50\%& 5.15\\
    \hline
  10 & 0.3537 & 0.2988 & 15.54\%& 0.2992& 15.41\%&4.09\\
    \hline
  20 & 0.3475 & 0.2912 & 16.20\%& 0.2928 & 15.74\%& 4.40\\
    \hline
  30 & 0.3454  & 0.2887 & 16.42\%& 0.2844 & 17.67\%& 3.45\\
    \hline
\end{tabular}
\end{center}
\end{table}

\begin{table}[H]
\caption{The IMSE and the reductions in the IMSE of $\hat{g}_n^{trap}$ using the uniform design, theoretical optimal design and estimated optimal design when $R(s,t)=\sigma^2 \rho^{|s^\lambda-t^\lambda|/\lambda}$ and $n=10.$ }\label{GMuTRAPo10}
\begin{center}
\begin{tabular}{|l|c|r|c|c|c|c|} 
  \hline
  $m$ & $Trap_{unif}$& $Trap_{opt}$& ${rIMSE}_{\lambda}$ & $Trap_{\widehat{opt}}$ & ${rIMSE}_{\widehat{\lambda}}$ & $\widehat{\lambda}$\\
  \hline
  5 & 0.1969  & 0.1771 & 10.06\%&  0.1822 & 7.50\%& 5.06\\
    \hline
  10 & 0.1674 & 0.1494 & 10.79\%& 0.1487 & 11.19\%& 3.91\\
    \hline
  20 & 0.1527 & 0.1355 & 11.26\%& 0.1305 & 14.54 \%& 3.21\\
    \hline
  30 & 0.1477 & 0.1309  & 11.43\%& 0.1346& 8.87\%&4.50 \\
    \hline
\end{tabular}
\end{center}
\end{table}

\begin{table}[H]
\caption{The IMSE and the reductions in the IMSE of $\hat{g}_n^{trap}$ using the uniform design, theoretical optimal design and estimated optimal design when $R(s,t)=\sigma^2 \rho^{|s^\lambda-t^\lambda|/\lambda}$ and $n=20.$ }\label{GMuTRAPo20}
\begin{center}
\begin{tabular}{|l|c|r|c|c|c|c|} 
  \hline
  $m$ & $Trap_{unif}$& $Trap_{opt}$&${rIMSE}_{\lambda}$ & $Trap_{\widehat{opt}}$ & ${rIMSE}_{\widehat{\lambda}}$ & $\widehat{\lambda}$\\
  \hline
  5 & 0.1699 & 0.1487 &  12.52\%& 0.1457 & 14.26\%& 4.35\\
  \hline
  10 & 0.1274 & 0.1096 & 12.14\%& 0.1106 & 11.34\%& 3.82\\
  \hline
  20 & 0.1022 & 0.0901 &  11.86\%& 0.0885& 13.39\%& 4.34\\
  \hline
  30 & 0.0947 & 0.0836 & 11.73\%& 0.0839 & 11.31\%& 3.90\\
  \hline
\end{tabular}
\end{center}
\end{table}

\begin{table}[H]
\caption{The IMSE and the reductions in the IMSE of $\hat{g}_n^{trap}$ using the uniform design, theoretical optimal design and estimated optimal design when $R(s,t)=\sigma^2 \rho^{|s^\lambda-t^\lambda|/\lambda}$ and $n=30.$ }\label{GMuTRAPo30}
\begin{center}
\begin{tabular}{|l|c|r|c|c|c|c|} 
  \hline
  $m$ & $Trap_{unif}$& $Trap_{opt}$& ${rIMSE}_{\lambda}$ & $Trap_{\widehat{opt}}$ & ${rIMSE}_{\widehat{\lambda}}$ & $\widehat{\lambda}$\\
  \hline
  5 &0.1682  & 0.1488 & 11.56\%& 0.1434 & 14.78\%& 4.46\\
  \hline
  10 & 0.1201 & 0.1056 & 12.09\%& 0.0973 &19.03\%& 4.86\\
  \hline
  20 & 0.0961 & 0.0840 & 12.57\%& 0.0861 & 10.4\%& 3.69\\
  \hline
  30 & 0.0881 & 0.0768 & 12.78\%& 0.7586& 13.88\%& 4.14\\
  \hline
\end{tabular}
\end{center}
\end{table}

\section{Proofs} 
\subsection{Proof of Lemma \ref{N_T=Onh}.}
For the sake of clarity, we omit the $n$ in $t_{i,n}$.  For $i=1,\cdots,n-1$ the Mean Value Theorem  (m.v.t) yields that there exists  $\eta_{i} \in ]t_i,t_{i+1}[$ such that,
\begin{align*}
t_{i+1}-t_i& =F^{-1}(\frac{i+1}{n})-F^{-1}(\frac{i}{n})= \frac{1}{nf(\eta_i)}.
\end{align*}
Since $\underset{0 \leq t \leq 1}{\inf} f(t) > 0$ then $t_{i+1}-t_i=O(\frac{1}{n})$. We shall now prove the second part of the Lemma. Since $T_n \cap [x-h,x+h] \ne \emptyset$, there exist $i_1, i_N$ indexes in $ \{1,\dots,n\}$ such that, \[N_{T_n} \leq i_N-i_1+1.\]
From the definition of the regular sequence we have for all $i = 1,...,n$, 
 \[t_i=F^{-1}\Big(\frac{i}{n}\Big) \;\;\;\text{thus}\;\;\; i = nF(t_i) .\]
Using this and the m.v.t we obtain for some $\epsilon_x \in ]t_{i_1},t_{i_N}[$,
 \begin{align*}
 N_{T_n}& \leq n\big(F(t_{i_N})-F(t_{i_1})\big)+1 = n(t_{i_N}-t_{i_1})f(\epsilon_{x})+1,
 \end{align*}
The boundedness of $f$ and the fact that $t_{i_N}-t_{i_1} \leq 2h$ yield,
\[N_{T_n} \leq (2~ \underset{0 \leq t \leq 1}{\sup}~ f(t))~ nh+1.\]
 This concludes the proof of the second part  of Lemma \ref{N_T=Onh} since   $ 1 \leq nh. $ $\Box$
\subsection{Proof of Proposition \ref{biaslemma}.}
For $h$ small enough and since $T_n \cap [x-h,x+h] \ne \emptyset$ we take $t_{x,1}<t_{x,2} < \cdots<  t_{x,N_{T_n}}$ the points of $T_n$ in $[x-h,x+h]$. 
Since $\mathbb{E}(\overline{Y}(t_{i}))=g(t_{i})$ for all $i=1,\cdots,n$ we have,
\begin{align*}
\mathbb{E} (\hat{g}^{\text{trap}}_n (x)) = \frac{1}{2n}\bigg\{  \sum_{k=1}^{N_{T_n}-1}\big( \frac{\varphi_{x,h}}{f}g\big) (t_{x,k})+ \big( \frac{\varphi_{x,h}}{f}g\big) (t_{x,k+1})\bigg\}.
\end{align*}
From the definition of the regular sequence of designs we have for $k=1,\dots,N_{T_n}-1$,
\begin{equation} \label{1/n=intf}
F(t_{x,k+1})-F(t_{x,k})=\frac{1}{n} \iff \int_{t_{x,k}}^{t_{x,k+1}}f(t)\;dt = \frac{1}{n}. 
\end{equation}Thus,
\begin{align*}
\mathbb{E} (\hat{g}^{\text{trap}}_n (x)) = \frac{1}{2} \sum_{k=1}^{N_{T_n}-1}\int_{t_{x,k}}^{t_{x,k+1}}\bigg\{ 
\big( \frac{\varphi_{x,h}}{f}g\big) (t_{x,k})+ \big( \frac{\varphi_{x,h}}{f}g\big) (t_{x,k+1})\bigg\}f(t)\;dt.
\end{align*}
Let, 
\begin{align*}
I_h(x)&=\int_{x-h}^{x+h}\varphi_{x,h}(t)g(t)\;dt \\
&= \sum_{k=1}^{N_{T_n}-1}\int_{t_{x,k}}^{t_{x,k+1}} \varphi_{x,h}(t)g(t)~dt + \int_{x-h}^{t_{x,1}} \varphi_{x,h}(t)g(t)~dt+\int_{t_{x,N_{T_n}}}^{x+h} \varphi_{x,h}(t)g(t)~dt,
\end{align*}
and write, 
\begin{align}
\mathbb{E} (\hat{g}^{\text{trap}}_n (x)) = \mathbb{E} (\hat{g}^{\text{trap}}_n (x))  -I_h(x)+I_h(x) \overset{\Delta}{=}\Delta_{x,h} +I_h(x). \label{E=}
\end{align}
We first control $\Delta_{x,h}$.  Let,
\begin{align}\label{deltabias}
\Delta_{x,h}  & = \Delta_{x,h}^1+\Delta_{x,h}^2,
\end{align}
where,
\begin{align*}
 \Delta_{x,h}^1 & = \frac{1}{2} \sum_{k=1}^{N_{T_n}-1}\int_{t_{x,k}}^{t_{x,k+1}}\bigg( \Big(\frac{\varphi_{x,h}}{f}g\Big)(t_{x,k}) f(t)-\varphi_{x,h}(t)g(t)\bigg)\;dt \\
 &~~~~~~~~~~~-\frac{1}{2}\int_{x-h}^{t_{x,1}}\varphi_{x,h}(t)g(t)~dt-\frac{1}{2}\int_{t_{x,N_{T_n}}}^{x+h}\varphi_{x,h}(t)g(t)~dt.\\
 \Delta_{x,h}^2 & = \frac{1}{2} \sum_{k=1}^{N_{T_n}-1}\int_{t_{x,k}}^{t_{x,k+1}}\bigg( \Big(\frac{\varphi_{x,h}}{f}g\Big)(t_{x,k+1}) f(t)-\varphi_{x,h}(t)g(t)\bigg)\;dt \\
 &~~~~~~~~~~~-\frac{1}{2}\int_{x-h}^{t_{x,1}}\varphi_{x,h}(t)g(t)~dt-\frac{1}{2}\int_{t_{x,N_{T_n}}}^{x+h}\varphi_{x,h}(t)g(t)~dt.
 \end{align*}For $t \in [x-h,t_{x,1}]$,  Taylor expansion of $\varphi_{x,h}$ around $(x-h)$ yields,
\begin{align}
\varphi_{x,h}(t) & = \varphi_{x,h}(x-h)+(t-(x-h))\varphi_{x,h}'(x-h)+\frac{1}{2}(t-(x-h))^2\varphi_{x,h}''(\theta_{x,h}),\label{varphiright}
\end{align}
for some $\theta_{x,h} \in ]x-h,t_{x,1}[$.  Recall that by definition of $\varphi_{x,h}$ we have,
\begin{equation}\label{phibounderies}
\underset{0 \leq t \leq 1}{\sup}~|\varphi_{x,h}^{(j)}(t)|\leq \frac{c_j}{h^{j+1}}~~\text{for}~j=0,1,2,
\end{equation}
for some appropriate constants $c_j$ where $j=0,1,2$. In addition, since  $\varphi_{x,h}$ is in $C^2$  and of support $[x-h,x+h]$ then, \begin{equation}\label{varphiis0}
\varphi_{x,h}(x-h)=\varphi_{x,h}(x+h)=\varphi_{x,h}'(x-h)=\varphi_{x,h}'(x+h)=0.
\end{equation}
Using \eqref{varphiis0} and \eqref{phibounderies} in \eqref{varphiright} and using Lemma \eqref{N_T=Onh} we obtain for  $t \in [x-h,t_{x,1}]$,
\begin{align}
\varphi_{x,h}(t) & = \frac{1}{2}(t-(x-h))^2\varphi_{x,h}''(\theta_{x,h}) = O\Big(\frac{1}{n^2h^3}\Big), \label{varphix-h}
\end{align}
 Likewise, for $t \in [t_{x,N_{T_n}},x+h]$ we have,
\begin{align}
\varphi_{x,h}(t) & = \frac{1}{2}(t-(x+h))^2\varphi_{x,h}''(\theta_{x,h}') = O\Big(\frac{1}{n^2h^3}\Big),\label{varphix+h}
\end{align}
where $\theta_{x,h}' \in ]t_{x,N_{T_n}},x+h[$. 
Hence, \begin{equation*}
\int_{x-h}^{t_{x,1}}\varphi_{x,h}(t)g(t)~dt=O\Big(\frac{1}{n^3h^3}\Big)~\text{and}~\int_{t_{x,N_{T_n}}}^{x+h}\varphi_{x,h}(t)g(t)~dt= O\Big(\frac{1}{n^3h^3}\Big).
\end{equation*}
Thus,
 \begin{align*}
 \Delta_{x,h}^1 & = \frac{1}{2} \sum_{k=1}^{N_{T_n}-1}\int_{t_{x,k}}^{t_{x,k+1}}\bigg( \Big(\frac{\varphi_{x,h}}{f}g\Big)(t_{x,k}) -\Big(\frac{\varphi_{x,h}}{f}g\Big)(t)\bigg)f(t)\;dt+O\Big(\frac{1}{n^3h^3}\Big),
\end{align*}
and,
\begin{align*}
 \Delta_{x,h}^2 & = \frac{1}{2} \sum_{k=1}^{N_{T_n}-1}\int_{t_{x,k}}^{t_{x,k+1}}\bigg( \Big(\frac{\varphi_{x,h}}{f}g\Big)(t_{x,k+1}) -\Big(\frac{\varphi_{x,h}}{f}g\Big)(t)\bigg)f(t)\;dt+O\Big(\frac{1}{n^3h^3}\Big).
\end{align*}
Recall that $\varphi_{x,h}$ is in $C^2$ and $f,g \in C^2([0,1])$, then for any $t \in ]t_{x,k},t_{x,k+1}[$ Taylor expansions of $\frac{\varphi_{x,h}}{f}g$ and $f$ around $t_{x,k}$ give,
\begin{align*}
 \Delta_{x,h}^1 & = \frac{1}{2} \sum_{k=1}^{N_{T_n}-1} \Big(\frac{\varphi_{x,h}}{f}g\Big)'(t_{x,k})f(t_{x,k}) \int_{t_{x,k}}^{t_{x,k+1}}(t_{x,k}-t)\;dt\\
&~~~~~- \frac{1}{2} \sum_{k=1}^{N_{T_n}-1} \Big(\frac{\varphi_{x,h}}{f}g\Big)'(t_{x,k})f'(t_{x,k}) \int_{t_{x,k}}^{t_{x,k+1}}(t-t_{x,k})^2\;dt\\
&~~~~~-\frac{1}{4} \sum_{k=1}^{N_{T_n}-1} \Big(\frac{\varphi_{x,h}}{f}g\Big)'(t_{x,k}) \int_{t_{x,k}}^{t_{x,k+1}}(t-t_{x,k})^3f''(\eta_{x,k})\;dt\\
 &~~~~~-\frac{1}{4}\sum_{k=1}^{N_{T_n}-1}f(t_{x,k})\int_{t_{x,k}}^{t_{x,k+1}}(t-t_{x,k})^2\Big(\frac{\varphi_{x,h}}{f}g\Big)''(\theta_{x,k})\;dt\\
 &~~~~~- \frac{1}{4}\sum_{k=1}^{N_{T_n}-1}f'(t_{x,k})\int_{t_{x,k}}^{t_{x,k+1}}(t-t_{x,k})^3\Big(\frac{\varphi_{x,h}}{f}g\Big)''(\theta_{x,k})\;dt\\
  &~~~~~-\frac{1}{8}\sum_{k=1}^{N_{T_n}-1}\int_{t_{x,k}}^{t_{x,k+1}}(t-t_{x,k})^4\Big(\frac{\varphi_{x,h}}{f}g\Big)''(\theta_{x,k})f''(\eta_{x,k})\;dt+O\Big(\frac{1}{n^3h^3}\Big),
\end{align*}
where $\theta_{x,k} $ and $\eta_{x,k} $ are in $ ]t_{x,k},t[$.
Recall that the functions $g^{(j)},f^{(j)}$ for $j=0,1,2$ are all bounded, then using\eqref{phibounderies} and  Lemma \ref{N_T=Onh} we get,
\begin{align}
&  \sum_{k=1}^{N_{T_n}-1} \Big(\frac{\varphi_{x,h}}{f}g\Big)'(t_{x,k}) \int_{t_{x,k}}^{t_{x,k+1}}(t-t_{x,k})^3f''(\eta_{x,k})\;dt =O\Big(\frac{1}{n^3h}\Big). \label{neg1}\\
& \sum_{k=1}^{N_{T_n}-1}f'(t_{x,k})\int_{t_{x,k}}^{t_{x,k+1}}(t-t_{x,k})^3\Big(\frac{\varphi_{x,h}}{f}g\Big)''(\theta_{x,k})\;dt=O\Big(\frac{1}{n^3h^2}\Big).\label{neg2}\\
& \sum_{k=1}^{N_{T_n}-1}\int_{t_{x,k}}^{t_{x,k+1}}(t-t_{x,k})^4\Big(\frac{\varphi_{x,h}}{f}g\Big)''(\theta_{x,k})f''(\eta_{x,k})\;dt=O\Big(\frac{1}{n^4h^2}\Big).\label{neg3}
\end{align}
Note that, since $\varphi_{x,h}''$, $g''$ and $f''$ are all lipschitz then,
\begin{align}
\Big(\frac{\varphi_{x,h}}{f}g\Big)''(\theta_{x,k}) & = \Big(\frac{\varphi_{x,h}}{f}g\Big)''(t_{x,k})+\bigg[ \Big(\frac{\varphi_{x,h}}{f}g\Big)''(\theta_{x,k})-\Big(\frac{\varphi_{x,h}}{f}g\Big)''(t_{x,k})\bigg]\nonumber\\
& = \Big(\frac{\varphi_{x,h}}{f}g\Big)''(t_{x,k})+O\Big(\frac{1}{nh^4}\Big). \label{phi''diff}
\end{align}
Injecting \eqref{neg1}, \eqref{neg2}, \eqref{neg3} and \eqref{phi''diff} in $\Delta_{x,h}^1$ we have, 
\begin{align*}
 \Delta_{x,h}^1 & = \frac{1}{2} \sum_{k=1}^{N_{T_n}-1} \Big(\frac{\varphi_{x,h}}{f}g\Big)'(t_{x,k})f(t_{x,k}) \int_{t_{x,k}}^{t_{x,k+1}}(t_{x,k}-t)\;dt\\
&~~~~~- \frac{1}{2} \sum_{k=1}^{N_{T_n}-1} \Big(\frac{\varphi_{x,h}}{f}g\Big)'(t_{x,k})f'(t_{x,k}) \int_{t_{x,k}}^{t_{x,k+1}}(t-t_{x,k})^2\;dt\\
 &~~~~~-\frac{1}{4}\sum_{k=1}^{N_{T_n}-1}\Big(\frac{\varphi_{x,h}}{f}g\Big)''(t_{x,k})f(t_{x,k})\int_{t_{x,k}}^{t_{x,k+1}}(t-t_{x,k})^2\;dt
+O\Big(\frac{1}{n^3h^3}\Big).
\end{align*}
Let $d_{x,k} = t_{x,k+1}-t_{x,k}$. We obtain by basic integration,
\begin{align}
 \Delta_{x,h}^1 & = -\frac{1}{4} \sum_{k=1}^{N_{T_n}-1} \Big(\frac{\varphi_{x,h}}{f}g\Big)'(t_{x,k})f(t_{x,k})d_{x,k}^2- \frac{1}{6} \sum_{k=1}^{N_{T_n}-1} \Big(\frac{\varphi_{x,h}}{f}g\Big)'(t_{x,k})f'(t_{x,k})d_{x,k}^3\nonumber \\
 &~~~~~-\frac{1}{12}\sum_{k=1}^{N_{T_n}-1}\Big(\frac{\varphi_{x,h}}{f}g\Big)''(t_{x,k})f(t_{x,k})d_{x,k}^3\;
+O\Big(\frac{1}{n^3h^3}\Big).\label{deltax1}
\end{align}
Similarly we verify that,
\begin{align}
 \Delta_{x,h}^2 & = \frac{1}{4} \sum_{k=1}^{N_{T_n}-1} \Big(\frac{\varphi_{x,h}}{f}g\Big)'(t_{x,k+1})f(t_{x,k+1})d_{x,k}^2- \frac{1}{6} \sum_{k=1}^{N_{T_n}-1} \Big(\frac{\varphi_{x,h}}{f}g\Big)'(t_{x,k+1})f'(t_{x,k+1})d_{x,k}^3 \nonumber  \\
 &~~~~~-\frac{1}{12}\sum_{k=1}^{N_{T_n}-1}\Big(\frac{\varphi_{x,h}}{f}g\Big)''(t_{x,k+1})f(t_{x,k+1})d_{x,k}^3\;
+O\Big(\frac{1}{n^3h^3}\Big),\label{deltax2}
\end{align}
Summing \eqref{deltax1} and \eqref{deltax2} gives,
\begin{align*}
&~~~~~~~~~~~~~~~~\Delta_{x,h}  = \Delta_{x,h}^1+\Delta_{x,h}^2\\
& = \frac{1}{4} \sum_{k=1}^{N_{T_n}-1} d_{x,k}^2  \bigg[ \Big(\frac{\varphi_{x,h}}{f}g\Big)'(t_{x,k+1})f(t_{x,k+1})-\Big(\frac{\varphi_{x,h}}{f}g\Big)'(t_{x,k})f(t_{x,k}) \bigg]\\
&~~~- \frac{1}{6} \sum_{k=1}^{N_{T_n}-1} d_{x,k}^3  \bigg[ \Big(\frac{\varphi_{x,h}}{f}g\Big)'(t_{x,k+1})f'(t_{x,k+1})+\Big(\frac{\varphi_{x,h}}{f}g\Big)'(t_{x,k})f'(t_{x,k})\bigg]\\
&~~~-\frac{1}{12}\sum_{k=1}^{N_{T_n}-1}d_{x,k}^3\bigg[\Big(\frac{\varphi_{x,h}}{f}g\Big)''(t_{x,k+1})f(t_{x,k+1})+\Big(\frac{\varphi_{x,h}}{f}g\Big)''(t_{x,k})f(t_{x,k})\bigg]+O\Big(\frac{1}{n^3h^3}\Big).
\end{align*}
Since $\varphi_{x,h}'$ is in $C^1$ and $g',f'\in C^1([0,1])$, Taylor expansion of $\Big(\frac{\varphi_{x,h}}{f}g\Big)'f$ around $t_{x,k}$ yields, \[\Big( \big(\frac{\varphi_{x,h}}{f}g\big)'f\Big)(t_{x,k+1}) = \Big( \big(\frac{\varphi_{x,h}}{f}g\big)'f\Big)(t_{x,k})+d_{x,k}\Big( \big(\frac{\varphi_{x,h}}{f}g\big)'f\Big)'(\nu_{x,k}),\]
where $\nu_{x,k} \in ]t_{x,k},t_{x,k+1}[$. We then have,
\begin{align*}
& \Delta_{x,h}  =  \frac{1}{4} \sum_{k=1}^{N_{T_n}-1} d_{x,k}^3\Big( \big(\frac{\varphi_{x,h}}{f}g\big)'f\Big)'(\nu_{x,k})  \\
&- \frac{1}{6} \sum_{k=1}^{N_{T_n}-1} d_{x,k}^3  \bigg[ \Big(\frac{\varphi_{x,h}}{f}g\Big)'(t_{x,k+1})f'(t_{x,k+1})+\Big(\frac{\varphi_{x,h}}{f}g\Big)'(t_{x,k})f'(t_{x,k})\bigg]\\
&-\frac{1}{12}\sum_{k=1}^{N_{T_n}-1}d_{x,k}^3\bigg[\Big(\frac{\varphi_{x,h}}{f}g\Big)''(t_{x,k+1})f(t_{x,k+1})+\Big(\frac{\varphi_{x,h}}{f}g\Big)''(t_{x,k})f(t_{x,k})\bigg]+O\Big(\frac{1}{n^3h^3}\Big).
\end{align*}
From the definition of the regular sequence of designs and using the m.v.t. we obtain for $k=1,\cdots,N_{T_n}-1,$
\begin{equation}\label{d=1/nf}
\int_{t_{x,k}}^{t_{x,k+1}}f(t)\;dt=\frac{1}{n}\iff d_{x,k}=\frac{1}{nf(t_{x,k}^*)}~\text{for some}~t_{x,k}^*\in]t_{x,k},t_{x,k+1}[.
\end{equation}
This equation yields,
\begin{align*}
\Delta_{x,h} & = \frac{1}{4n^2} \sum_{k=1}^{N_{T_n}-1} d_{x,k}  \frac{1}{f^2(t_{x,k}^*)}\Big( \big(\frac{\varphi_{x,h}}{f}g\big)'f\Big)'(\nu_{x,k})\\
&~~~- \frac{1}{6n^2} \sum_{k=1}^{N_{T_n}-1} d_{x,k}  \bigg[ \Big(\frac{\varphi_{x,h}}{f}g\Big)'(t_{x,k+1})\frac{f'(t_{x,k+1})}{f^2(t_{x,k}^*)}+\Big(\frac{\varphi_{x,h}}{f}g\Big)'(t_{x,k})\frac{f'(t_{x,k})}{f^2(t_{x,k}^*)}\bigg]\\
&~~~-\frac{1}{12n^2}\sum_{k=1}^{N_{T_n}-1}d_{x,k}\bigg[\Big(\frac{\varphi_{x,h}}{f}g\Big)''(t_{x,k+1})\frac{f(t_{x,k+1})}{f^2(t_{x,k}^*)}+\Big(\frac{\varphi_{x,h}}{f}g\Big)''(t_{x,k})\frac{f(t_{x,k})}{f^2(t_{x,k}^*)}\bigg]\\
&~~~+O\Big(\frac{1}{n^3h^3}\Big).
\end{align*}
Using the Riemann integrability of $\varphi_{x,h}^{(j)},f^{(j)}$ and $g^{(j)}$ for $j=0,1,2$ and applying  Lemma  \ref{sumtointegral} in the Appendix with $u(t)=\frac{1}{f^2(t)}$ and $v(t)=\Big( \big(\frac{\varphi_{x,h}}{f}g\big)'f\Big)'(t)$ we obtain,
\begin{align*}
 \sum_{k=1}^{N_{T_n}-1} d_{x,k}  \frac{1}{f^2(t_{x,k}^*)}\Big( \big(\frac{\varphi_{x,h}}{f}g\big)'f\Big)'(\nu_{x,k}) & =  \int_{x-h}^{x+h}  \frac{1}{f^2(t)}\Big( \big(\frac{\varphi_{x,h}}{f}g\big)'f\Big)'(t)~dt+O\Big(\frac{1}{nh^3}\Big).
\end{align*}
Similarly, taking $u(t)=\Big(\frac{\varphi_{x,h}}{f}g\Big)'(t)$ and $v(t)=\frac{f'(t)}{f^2(t)}$ in Lemma \ref{sumtointegral} we obtain,
\begin{align*}
 \sum_{k=1}^{N_{T_n}-1} d_{x,k}  \Big(\frac{\varphi_{x,h}}{f}g\Big)'(t_{x,k+1})\frac{f'(t_{x,k+1})}{f^2(t_{x,k}^*)} & = \int_{x-h}^{x+h} \Big(\frac{\varphi_{x,h}}{f}g\Big)'(t)\frac{f'(t)}{f^2(t)}~dt+O\Big(\frac{1}{nh^3}\Big).
\end{align*}
Again taking $u(t)=\Big(\frac{\varphi_{x,h}}{f}g\Big)''(t)$ and $v(t)=\frac{f(t)}{f^2(t)}$ we obtain,
\begin{align*}
\sum_{k=1}^{N_{T_n}-1}d_{x,k}\Big(\frac{\varphi_{x,h}}{f}g\Big)''(t_{x,k+1})\frac{f(t_{x,k+1})}{f^2(t_{x,k}^*)}=\int_{x-h}^{x+h} \Big(\frac{\varphi_{x,h}}{f}g\Big)''(t)\frac{1}{f(t)}~dt+O\Big(\frac{1}{nh^3}\Big).
\end{align*}
Hence,
\begin{align*}
\Delta_{x,h} & = \frac{1}{4n^2} \int_{x-h}^{x+h}  \frac{1}{f^2(t)}\Big( \big(\frac{\varphi_{x,h}}{f}g\big)'f\Big)'(t)~dt- \frac{1}{3n^2} \int_{x-h}^{x+h} \Big(\frac{\varphi_{x,h}}{f}g\Big)'(t)\frac{f'(t)}{f^2(t)}~dt\\
&~~~-\frac{1}{6n^2}\int_{x-h}^{x+h} \Big(\frac{\varphi_{x,h}}{f}g\Big)''(t)\frac{1}{f(t)}~dt+O\Big(\frac{1}{n^3h^3}\Big).
\end{align*}
Simple derivations yield,
\begin{align*}
& \Delta_{x,h}  = \frac{1}{4n^2} \int_{x-h}^{x+h}  \Big(\frac{\varphi_{x,h}}{f}g\Big)''(t)\frac{1}{f(t)}~dt+\frac{1}{4n^2} \int_{x-h}^{x+h}\Big(\frac{\varphi_{x,h}}{f}g\Big)'(t)\frac{f'(t)}{f^2(t)}~dt \\
&- \frac{1}{3n^2} \int_{x-h}^{x+h} \Big(\frac{\varphi_{x,h}}{f}g\Big)'(t)\frac{f'(t)}{f^2(t)}~dt-\frac{1}{6n^2}\int_{x-h}^{x+h} \Big(\frac{\varphi_{x,h}}{f}g\Big)''(t)\frac{1}{f(t)}~dt+O\Big(\frac{1}{n^3h^3}\Big)\\
& = \frac{1}{12n^2} \int_{x-h}^{x+h}  \Big(\frac{\varphi_{x,h}}{f}g\Big)''(t)\frac{1}{f(t)}~dt- \frac{1}{12n^2} \int_{x-h}^{x+h} \Big(\frac{\varphi_{x,h}}{f}g\Big)'(t)\frac{f'(t)}{f^2(t)}~dt+O\Big(\frac{1}{n^3h^3}\Big)\\
& = \frac{1}{12n^2} \int_{x-h}^{x+h}  \bigg( \Big(\frac{\varphi_{x,h}}{f}g\Big)'\frac{1}{f}\bigg)'(t)~dt+O\Big(\frac{1}{n^3h^3}\Big).
\end{align*}
Finally,
\begin{align*}
\Delta_{x,h} & =\frac{1}{12n^2}\bigg( \Big(\frac{\varphi_{x,h}}{f}g\Big)'\frac{1}{f}\bigg)(x+h) - \bigg( \Big(\frac{\varphi_{x,h}}{f}g\Big)'\frac{1}{f}\bigg)(x-h)+O\Big(\frac{1}{n^3h^3}\Big).
\end{align*}
The last equation together with \eqref{varphiis0} yield, 
\begin{equation}
\Delta_{x,h}=O\Big(\frac{1}{n^3h^3}\Big). \label{Delta=}
\end{equation}
The control of $I_h(x)$ is classical and it can be seen from Gasser and Müller  \cite{Gasser and Muller 1984} that,
\begin{align}
 I_h(x) & = g(x)+ \frac{1}{2}h^2 g''(x)\int_{-1}^{1}t^2K(t) ~dt + o (h^2). \label{lhx=}
\end{align}
Finally, collecting \eqref{E=}, \eqref{Delta=} and \eqref{lhx=} gives, \[\Bias(\hat{g}^{trap}_n (x) )= \frac{1}{2} h^2 g''(x)B + o (h^2) + O\Big(\frac{1}{n^3h^3}\Big),\]
where $B = \int_{-1}^{1} t^2 K(t)~dt $. 
This concludes the proof of Proposition \ref{biaslemma}. $\Box$
\subsection{Proof of Proposition \ref{variancelemma1}.} The greatest lines of this proof are based on the work of Belouni and Benhenni \cite{Belouni}. For $h$ small enough and since $T_n \cap [x-h,x+h] \ne \emptyset$ we have, \[0 \leq t_1 < \cdots < x-h \leq t_{x,1}< \dots < t_{x,N_{T_n}} \leq x+h <\dots <t_n\leq1. \] 
Let,\[\Phi(t,s) = \big( \frac{\varphi_{x,h}}{f}\big )(t) R(t,s)\big( \frac{\varphi_{x,h}}{f}\big )(s),\]
and,
\begin{equation} \label{sigma_x,h}
\sigma_{x,h}^2 =  \int_{x-h}^{x+h}\int_{x-h}^{x+h} \varphi_{x,h}(t)R(t,s)\varphi_{x,h}(s)~ds~dt.
\end{equation}
On the one hand,
\begin{align*}
\Var(\hat{g}^{trap}_n (x)) &= \frac{1}{4mn^2} \sum_{i=1}^{N_{T_n}-1}\sum_{j=1}^{N_{T_n}-1} \bigg\{ \Phi(t_{x,i},t_{x,j})+\Phi(t_{x,i},t_{x,j+1})+\Phi(t_{x,i+1},t_{x,j})\\
&~~~~~~~~~~~~~+\Phi(t_{x,i+1},t_{x,j+1}) \bigg\}
\end{align*}
Using \eqref{1/n=intf} one can write,
\begin{align*}
\Var(\hat{g}^{trap}_n (x)) = \frac{1}{4m} \sum_{i=1}^{N_{T_n}-1}\sum_{j=1}^{N_{T_n}-1} \int_{t_{x,i}}^{t_{x,i+1}}\int_{t_{x,j}}^{t_{x,j+1}}\bigg\{ \Phi(t_{x,i},t_{x,j})+\Phi(t_{x,i},t_{x,j+1})\\
~~~~~~~~~~~~~~~~~~~~+\Phi(t_{x,i+1},t_{x,j})+\Phi(t_{x,i+1},t_{x,j+1}) \bigg\}f(s)\;f(t)\;ds\;dt.
\end{align*}
On the other hand we have,
{\small
\begin{align*}
\sigma_{x,h}^2 & =  \sum_{i=1}^{N_{T_n}-1}\sum_{j=1}^{N_{T_n}-1}  \int_{t_{x,i}}^{t_{x,i+1}}\int_{t_{x,j}}^{t_{x,j+1}}\Phi(t,s) f(t)\;f(s)\;ds\;dt\\
& + 2\int_{x-h}^{t_{x,1}}\int_{t_{x,N_{T_n}}}^{x+h}\Phi(t,s) f(t)\;f(s)\;ds\;dt + \int_{t_{x,N_{T_n}}}^{x+h}\int_{t_{x,N_{T_n}}}^{x+h}\Phi(t,s)f(t)\;f(s)\;ds\;dt\\
&+\int_{x-h}^{t_{x,1}}\int_{x-h}^{t_{x,1}}\Phi(t,s) f(t)\;f(s)\;ds\;dt+2\sum_{j=1}^{N_{T_n}-1} \int_{x-h}^{t_{x,1}}\int_{t_{x,j}}^{t_{x,j+1}}\Phi(t,s) f(t)\;f(s)\;ds\;dt\\
&+2\sum_{j=1}^{N_{T_n}-1} \int_{t_{x,N_{T_n}}}^{x+h}\int_{t_{x,j}}^{t_{x,j+1}}\Phi(t,s) f(t)\;f(s)\;ds\;dt.
\end{align*}}
Recall that Lemma \ref{N_T=Onh} yields $N_{T_n}=O(nh)$ and $\underset{1 \leq i \leq n}{\sup} d_{x,i} = O(\frac{1}{n})$. Using \eqref{varphix-h} and \eqref{varphix+h} we have, \begin{equation}\label{phiclosetoboundaries}\underset{(x-h) \leq t \leq t_{x,1}}{\sup} |\varphi_{x,h} (t)| = O\Big(\frac{1}{n^2h^3}\Big)~~~\text{and}~~\underset{ t_{x,N_{T_n}} \leq t \leq (x+h)}{\sup} |\varphi_{x,h} (t) |= O\Big(\frac{1}{n^2h^3}\Big).\end{equation}
Since $f$ and $ R$ are bounded, using \eqref{phibounderies} and \eqref{phiclosetoboundaries} we obtain,
\begin{align*}
& \int_{x-h}^{t_{x,1}}\int_{t_{x,N_{T_n}}}^{x+h}\Phi(t,s) f(t)\;f(s)\;ds\;dt  = O\Big(\frac{1}{n^6h^6}\Big),\\
& \int_{t_{x,N_{T_n}}}^{x+h}\int_{t_{x,N_{T_n}}}^{x+h}\Phi(t,s) f(t)\;f(s)\;ds\;dt= O\Big(\frac{1}{n^6h^6}\Big),\\
& \int_{x-h}^{t_{x,1}}\int_{x-h}^{t_{x,1}}\Phi(t,s) f(t)\;f(s)\;ds\;dt = O\Big(\frac{1}{n^6h^6}\Big),\\
 & \sum_{j=1}^{N_{T_n}-1} \int_{x-h}^{t_{x,1}}\int_{t_{x,j}}^{t_{x,j+1}}\Phi(t,s) f(t)\;f(s)\;ds\;dt = O\Big(\frac{1}{n^3h^3}\Big),\\
& \sum_{j=1}^{N_{T_n}-1} \int_{t_{x,N_{T_n}}}^{x+h}\int_{t_{x,j}}^{t_{x,j+1}}\Phi(t,s) f(t)\;f(s)\;ds\;dt= O\Big(\frac{1}{n^3h^3}\Big).
\end{align*}
Thus,\[\sigma_{x,h}^2 = \sum_{i=1}^{N_{T_n}-1}\sum_{j=1}^{N_{T_n}-1} \int_{t_{x,i}}^{t_{x,i+1}}\int_{t_{x,j}}^{t_{x,j+1}}\Phi(t,s) f(t)\;f(s)\;ds\;dt +O\Big(\frac{1}{n^3h^3}\Big). \]
We shall control the residual variance $ \Var(\hat{g}^{trap}_n (x))-\frac{\sigma_{x,h}^2}{m} $. For this, let, 
\begin{equation}\label{Nij} 
N_{i,j} (t,s) = \Phi(t_{x,i},t_{x,j})+\Phi(t_{x,i+1},t_{x,j})+\Phi(t_{x,i},t_{x,j+1})+\Phi(t_{x,i+1},t_{x,j+1})- 4 \Phi(t,s),
\end{equation}
and put,
\begin{align}
I_{i,j} & = \frac{1}{4m}\int_{t_{x,i}}^{t_{x,i+1}}\int_{t_{x,j}}^{t_{x,j+1}} N_{i,j} (t,s)f(t)\;f(s)\;ds\;dt. \label{Iijdef}
\end{align}The residual variance can then be written as follows, 
\begin{align}\label{var-sigma2}
\Var(\hat{g}^{trap}_n (x))-\frac{\sigma_{x,h}^2}{m} & = \sum_{i=1}^{N_{T_n}-1}  I_{i,i} + \overset{{N_{T_n}-1}}{\underset{i \ne j=1}{\sum\sum}  } I_{i,j}+O\Big(\frac{1}{mn^3h^3}\Big),
\end{align} 
Starting with the diagonal terms $I_{i,i}$. Since for any $s,t \in [0,1]$, we have $N_{i,i}(s,t)=N_{i,i}(t,s)$, then we can write,
\begin{align}
I_{i,i} & =\frac{1}{2m} \int_{t_{x,i}}^{t_{x,i+1}}\int_{t_{x,i}}^{t}N_{i,i} (t,s) f(t)\;f(s)\;ds\;dt. \label{1Iii}
\end{align}
Because of Assumption $(B)$, $N_{i,i}$ has left and right first order derivatives on the diagonal on $[0,1]^2$. For any $s,t$ such that ($t_{x,i} < s \leq t < t_{x,i+1}$), Taylor expansion of $\Phi$  around $(t_{x,i},t_{x,i})$ gives,
\begin{align*}
 \Phi(t,s)  &= \Phi(t,t_{x,i})+(s-t_{x,i})\Phi^{(0,1)}(t,t_{x,i})+\frac{1}{2}(s-t_{x,i})^2\Phi^{(0,2)}(t,\eta_{s,i}^{(1)}) \nonumber\\
& = \Phi(t_{x,i},t_{x,i})+(t-t_{x,i})\Phi^{(1,0)}(\epsilon_{t,i}^{(1)},t_{x,i})+(s-t_{x,i})\Phi^{(0,1)}(\epsilon_{i},t_{x,i})\nonumber\\
&~~~~+(s-t_{x,i})(t-\epsilon_{i})\Phi^{(1,1)}(\epsilon_{t,i}^{(2)},t_{x,i})+\frac{1}{2}(s-t_{x,i})^2\Phi^{(0,2)}(t,\eta_{s,i}^{(1)}),
\end{align*}
for some $\epsilon_{i} \in ]t_{x,i},t_{x,i+1}[$, some $\epsilon_{t,i}^{(1)}$ in  $]t_{x,i},t[$, some $\epsilon_{t,i}^{(2)}$ between $t$ and $\epsilon_{i}$ and some $\eta_{t,i}^{(1)}$ in $]t_{x,i},s [$. We have,
\begin{align*}
 \Phi(t,s) &= \Phi(t_{x,i},t_{x,i})+(t-t_{x,i})\Phi^{(1,0)}(\epsilon_{i},t_{x,i})+(s-t_{x,i})\Phi^{(0,1)}(\epsilon_{i},t_{x,i})\nonumber\\
&+(t-t_{x,i})\Big( \Phi^{(1,0)}(\epsilon_{t,i}^{(1)},t_{x,i})-\Phi^{(1,0)}(\epsilon_{i},t_{x,i})\Big)\\
&+(s-t_{x,i})(t-\epsilon_{i})\Phi^{(1,1)}(\epsilon_{t,i}^{(2)},t_{x,i})+\frac{1}{2}(s-t_{x,i})^2\Phi^{(0,2)}(t,\eta_{s,i}^{(1)}).
\end{align*}
For $l$ and $l'$ integers such that $l+l' \leq 2$, Assumption $(C)$ yields,
\begin{equation}\label{O20O11}
\underset{s \ne t }{\sup}~|\Phi^{(l,l')}(t,s)|=O\Big(\frac{1}{h^{l+l'+2}}\Big).
\end{equation}
In addition, since $\varphi_{x,h}, \varphi_{x,h}', \frac{1}{f}, R$ and $R(\cdot,t_{x,i})$  are all continuous on $]t_{x,i},t_{x,i+1}[$, then fo $s\ne t$ in $]t_{x,i},t_{x,i+1}[$ we have, 
 \begin{align*}
&\Big| \Phi^{(1,0)}(s,t_{x,i})-\Phi^{(1,0)}(t,t_{x,i})\Big|=\Big|\frac{\varphi_{x,h}}{f}(t_{x,i})\Big|\Big|R(s,t_{x,i})\bigg( \frac{\varphi_{x,h}'}{f}(s)-\frac{\varphi_{x,h}'}{f}(t)\bigg) \\
&+R^{(1,0)}(s,t_{x,i})\bigg( \frac{\varphi_{x,h}}{f}(s)-\frac{\varphi_{x,h}}{f}(t)\bigg)+\frac{\varphi_{x,h}}{f}(t)\Big(R^{(1,0)}(s,t_{x,i})-R^{(1,0)}(t,t_{x,i}) \Big)\\
&+\frac{\varphi_{x,h}'}{f}(t)\Big(R(s,t_{x,i})-R(t,t_{x,i}) \Big)\Big|=O\Big ( \frac{1}{nh^4}\Big).
\end{align*}
Finally, using this equation together with Lemma \ref{N_T=Onh} we obtain,
\begin{align}
 \Phi(t,s) &= \Phi(t_{x,i},t_{x,i})+(t-t_{x,i})\Phi^{(1,0)}(\epsilon_{i},t_{x,i})+(s-t_{x,i})\Phi^{(0,1)}(\epsilon_{i},t_{x,i})+O\Big ( \frac{1}{n^2h^4}\Big).\label{Kst}
\end{align}
Similarly we verify that,
\begin{align}
&\Phi(t_{x,i+1},t_{x,i+1}) = \Phi(t_{x,i},t_{x,i})+d_{x,i} \Phi^{(1,0)}(\epsilon_{i},t_{x,i})+d_{x,i} \Phi^{(0,1)}(\epsilon_{i},t_{x,i})+O\Big ( \frac{1}{n^2h^4}\Big),\label{Ki1i1}
\end{align}
and that,
\begin{align}
\Phi(t_{x,i+1},t_{x,i}) & = \Phi(t_{x,i},t_{x,i})+d_{x,i} \Phi^{(1,0)}(\epsilon_{i},t_{x,i})+O\Big ( \frac{1}{n^2h^4}\Big).\label{Ki1i}
\end{align}
Inserting \eqref{Kst}, \eqref{Ki1i1} and \eqref{Ki1i} in \eqref{Nij} for $i=j$ and using \eqref{O20O11} and Lemma \ref{N_T=Onh}, we obtain, 
\begin{align*}
N_{i,i} (t,s) & = 3 d_{x,i} \Phi^{(1,0)}(\epsilon_{i},t_{x,i})- 4 (t-t_{x,i}) \Phi^{(1,0)}(\epsilon_{i},t_{x,i})\\
&~~~~+d_{x,i} \Phi^{(0,1)}(\epsilon_{i},t_{x,i}) \big) - 4 (s-t_{x,i}) \Phi^{(0,1)}(\epsilon_{i},t_{x,i})+O\Big(\frac{1}{n^2h^4}\Big).
\end{align*}
Replacing this expression in \eqref{1Iii}, and using the boundedness of $f$ and Lemma \ref{N_T=Onh}, we obtain,
\begin{align}
I_{i,i} & = \frac{1}{2m}\Bigg( d_{x,i} \Big(3\Phi^{(1,0)}(\epsilon_{i},t_{x,i}) +\Phi^{(0,1)}(\epsilon_{i},t_{x,i})\Big) \int_{t_{x,i}}^{t_{x,i+1}}\int_{t_{x,i}}^{t}f(t)\;f(s)\;ds\;dt \nonumber\\
& ~~ -4\Phi^{(1,0)}(\epsilon_{i},t_{x,i}) \int_{t_{x,i}}^{t_{x,i+1}}\int_{t_{x,i}}^{t}(t-t_{x,i}) f(t)\;f(s)\;ds\;dt  \nonumber\\
& ~~ -4\Phi^{(0,1)}(\epsilon_{i},t_{x,i}) \int_{t_{x,i}}^{t_{x,i+1}}\int_{t_{x,i}}^{t}(s-t_{x,i}) f(t)\;f(s)\;ds\;dt \Bigg)+O\Big(\frac{1}{mn^4h^4}\Big). \label{Iiiii}
\end{align}
Recall that $f$ is  in $C^2([0,1])$ and that $d_{x,i}=O(\frac{1}{n})$ from Lemma  \ref{N_T=Onh}. It can easily be verified that for any integers  $l$ and $l'$:
\begin{align*}
& \int_{t_{x,i}}^{t_{x,i+1}} \int_{t_{x,i}}^{t} (s-t_{x,i})^{l'} (s-t_{x,i})^l f(t)\;f(s)\;ds\;dt = \frac{f^2(t_{x,i})~d_{x,i}^{(l+l'+2)}}{(l'+1)(l+l'+2)}+O\Big(\frac{1}{n^{l+l'+3}} \Big).
\end{align*}
Using this last Equation together with \eqref{O20O11}  in \eqref{Iiiii} above, and \eqref{d=1/nf} we obtain,
\begin{align*}
I_{i,i} & =  \frac{1}{12m}\Big(\Phi^{(1,0)}(\epsilon_{i},t_{x,i})-\Phi^{(0,1)}(\epsilon_{i},t_{x,i})\Big)f^2(t_{x,i})d_{x,i}^3  +O\Big( \frac{1}{mn^4h^4}\Big)\\
 & =  \frac{1}{12mn^2}\Big(\Phi^{(1,0)}(\epsilon_{i},t_{x,i})-\Phi^{(0,1)}(\epsilon_{i},t_{x,i})\Big)\frac{f^2(t_{x,i})}{f^2(t_{x,i}^*)}d_{x,i}  +O\Big( \frac{1}{mn^4h^4}\Big).
\end{align*}
Finally using Lemma  \ref{N_T=Onh}, the integrability of $\varphi_{x,h},\varphi_{x,h}', f, f'$ and $R^{(0,1)}(.,t)$ and applying Lemma \ref{sumtointegral} in the Appendix, we obtain,
{\small
\begin{align}
\sum_{i=1}^{N_{T_n}-1} I_{i,i} & = \frac{1}{12mn^2} \sum_{i=1}^{N_{T_n}-1}  \Big(\Phi^{(1,0)}(\epsilon_{i},t_{x,i})-\Phi^{(0,1)}(\epsilon_{i},t_{x,i})\Big)\frac{f^2(t_{x,i})}{f^2(t_{x,i}^*)}d_{x,i}+O\Big( \frac{1}{mn^3h^3}\Big)  \nonumber\\
& =  \frac{1}{12mn^2} \int_{x-h}^{x+h}\Big( \Phi^{(1,0)}(t^+,t) -\Phi^{(0,1)}(t^+,t)\Big) \;dt+ O\Big(\frac{1}{mn^3hÂ^3}\Big).
\end{align}}
Since $\Phi^{(0,1)}(t^+,t)=\Phi^{(0,1)}(t,t^-)=\Phi^{(1,0)}(t^-,t)$, then,
\begin{align}
\sum_{i=1}^{N_{T_n}-1} I_{i,i} 
& =  -\frac{1}{12mn^2} \int_{x-h}^{x+h}\Big(\Phi^{(1,0)}(t^-,t)- \Phi^{(1,0)}(t^+,t) \Big) \;dt+ O\Big(\frac{1}{mn^3hÂ^3}\Big).\label{diagtermtrap}
\end{align}
Now, it remains to handle the off diagonal term. 
Assumption $(B)$ yields that $N_{i,j}$ for $i \ne j$ is twice differentiable off the diagonal on $[0,1]^2$. Taylor expansion of $N_{i,j}$ around $(t_{x,i},t_{x,j})$ for $i \ne j$ up to order 2 gives,
\begin{align}
\Phi(t,s) & = \Phi(t_{x,i},t_{x,j})+(t-t_{x,i})\Phi^{(1,0)}(t_{x,i},t_{x,j})+(s-t_{x,j})\Phi^{(0,1)}(t_{x,i},t_{x,j})\nonumber\\
&~~~~~~~ +\frac{1}{2}(t-t_{x,i})^2\Phi^{(2,0)}(\epsilon_{x,i}^{(1)},t_{x,j})+\frac{1}{2}(s-t_{x,j})^2\Phi^{(0,2)}(t_{x,i},\eta_{x,j}^{(1)})\nonumber\\
&~~~~~~~  + (t-t_{x,i})(s-t_{x,j})\Phi^{(1,1)}(\epsilon_{x,i}^{(1)},\eta_{x,j}^{(1)}),
\label{kts'}
\end{align}
for some $\epsilon_{x,i}^{(1)}$ between  $t_{x,i}$ and $t$ and some $\eta_{x,j}^{(1)}$ between  $t_{x,j}$ and $s$. Taking $t=t_{x,i+1}$ and $s=t_{x,j}$ in \eqref{kts'}, we obtain,
\begin{align}
\Phi(t_{x,i+1},t_{x,j}) & = \Phi(t_{x,i},t_{x,j})+d_{x,i} \Phi^{(1,0)}(t_{x,i},t_{x,j})+\frac{1}{2}d_{x,i}^2\Phi^{(2,0)}(\epsilon_{x,i}^{(2)},t_{x,j}), \label{Kti1tj'}
\end{align}
for some $\epsilon_{x,i}^{(2)}$ in $]t_{x,i},t_{x,i+1}[$. Taking $t=t_{x,i}$ and $s=t_{x,j+1}$ in \eqref{kts'}, we obtain,
\begin{align}
\Phi(t_{x,i},t_{x,j+1}) & = \Phi(t_{x,i},t_{x,j})+d_{x,j} \Phi^{(0,1)}(t_{x,i},t_{x,j})+\frac{1}{2}d_{x,j}^2\Phi^{(0,2)}(t_{x,i},\eta_{x,j}^{(2)}), \label{Ktitj1'}
\end{align}
for some $\eta_{x,j}^{(2)}$ in $]t_{x,j},t_{x,j+1}[$. Taking $t=t_{x,i+1}$ and $s=t_{x,j+1}$ in \eqref{kts'}, we obtain,
\begin{align}
\Phi(t_{x,i+1},t_{x,j+1}) & = \Phi(t_{x,i},t_{x,j})+d_{x,i} \Phi^{(1,0)}(t_{x,i},t_{x,j})+d_{x,j} \Phi^{(0,1)}(t_{x,i},t_{x,j})\nonumber\\
&~~~+\frac{1}{2}d_{x,i}^2\Phi^{(2,0)}(\epsilon_{x,i}^{(3)},t_{x,j})  +\frac{1}{2}d_{x,j}^2\Phi^{(0,2)}(t_{x,i},\eta_{x,j}^{(3)}) \nonumber\\
&~~~+ d_{x,i} d_{x,j}\Phi^{(1,1)}(\epsilon_{x,i}^{(3)},\eta_{x,j}^{(3)}),\label{kti1tj1'}
\end{align}
We obtain by inserting \eqref{kts'}, \eqref{Kti1tj'}, \eqref{Ktitj1'}  and \eqref{kti1tj1'} in \eqref{Nij},
\begin{align*}
N_{i,j} (t,s) & =  \Phi^{(1,0)}(t_{x,i},t_{x,j}) \big(2d_{x,i} -4(t-t_{x,i})\big)+\Phi^{(0,1)}(t_{x,i},t_{x,j}) \big (2d_{x,j} -4(s-t_{x,j})\big)\\
&~~~~ + \frac{1}{2}d_{x,i}^2 \big(\Phi^{(2,0)}(\epsilon_{x,i}^{(2)},t_{x,j})+ \Phi^{(2,0)}(\epsilon_{x,i}^{(3)},t_{x,j})\big)-2(t-t_{x,i})^2\Phi^{(2,0)}(\epsilon_{x,i}^{(1)},t_{x,j})\\
&~~~~ + \frac{1}{2}d_{x,j}^2\big( \Phi^{(0,2)}(t_{x,i},\eta_{x,j}^{(2)})+\Phi^{(0,2)}(t_{x,i},\eta_{x,j}^{(3)})\big) -2(s-t_{x,j})^2\Phi^{(0,2)}(t_{x,i},\eta_{x,j}^{(1)})\\
&~~~~ + d_{x,i} d_{x,j}~\Phi^{(1,1)}(\epsilon_{x,i}^{(3)},\eta_{x,j}^{(3)})-4(t-t_{x,i})(s-t_{x,j})\Phi^{(1,1)}(\epsilon_{x,i}^{(1)},\eta_{x,j}^{(1)}).
\end{align*}
We obtain inserting the last equation in \eqref{Iijdef},
\begin{align}
I_{i,j}  & = \frac{1}{4m} \sum_{l=1}^{5} I_{i,j}^{(l)}, \label{Iij}
\end{align}
where,
\begin{align*}
I_{i,j}^{(1)} & = \Phi^{(1,0)}(t_{x,i},t_{x,j}) \bigg(2d_{x,i}\int_{t_{x,i}}^{t_{x,i+1}}\int_{t_{x,j}}^{t_{x,j+1}}f(t)f(s)dtds \\
&~~~~-4\int_{t_{x,i}}^{t_{x,i+1}}\int_{t_{x,j}}^{t_{x,j+1}}(t-t_{x,i})f(t)f(s)dtds\bigg).\\
I_{i,j}^{(2)} & = \Phi^{(0,1)}(t_{x,i},t_{x,j}) \bigg(2d_{x,j}\int_{t_{x,i}}^{t_{x,i+1}}\int_{t_{x,j}}^{t_{x,j+1}}f(t)f(s)dtds \\
&~~~~-4\int_{t_{x,i}}^{t_{x,i+1}}\int_{t_{x,j}}^{t_{x,j+1}}(s-t_{x,j})f(t)f(s)dtds\bigg).\\
I_{i,j}^{(3)} & =\frac{1}{2}d_{x,i}^2 \int_{t_{x,i}}^{t_{x,i+1}}\int_{t_{x,j}}^{t_{x,j+1}}\big(\Phi^{(2,0)}(\epsilon_{x,i}^{(2)},t_{x,j})+ \Phi^{(2,0)}(\epsilon_{x,i}^{(3)},t_{x,j})\big)f(t)f(s)dtds\\
&~~~~~~-2\int_{t_{x,i}}^{t_{x,i+1}}\int_{t_{x,j}}^{t_{x,j+1}}(t-t_{x,i})^2\Phi^{(2,0)}(\epsilon_{x,i}^{(1)},t_{x,j})f(t)f(s)dtds.\\
I_{i,j}^{(4)} & =\frac{1}{2}d_{x,j}^2 \int_{t_{x,i}}^{t_{x,i+1}}\int_{t_{x,j}}^{t_{x,j+1}}\big(\Phi^{(0,2)}(t_{x,i},\eta_{x,j}^{(2)})+\Phi^{(0,2)}(t_{x,i},\eta_{x,j}^{(3)})\big)f(t)f(s)dtds\\
&~~~~~~-2\int_{t_{x,i}}^{t_{x,i+1}}\int_{t_{x,j}}^{t_{x,j+1}}(s-t_{x,j})^2\Phi^{(0,2)}(t_{x,i},\eta_{x,j}^{(1)})f(t)f(s)dtds.\\
I_{i,j}^{(5)} & = d_{x,i} d_{x,j}\int_{t_{x,i}}^{t_{x,i+1}}\int_{t_{x,j}}^{t_{x,j+1}}\Phi^{(1,1)}(\epsilon_{x,i}^{(3)},\eta_{x,j}^{(3)})f(t)f(s)dtds\\
&~~~~~~-4\int_{t_{x,i}}^{t_{x,i+1}}\int_{t_{x,j}}^{t_{x,j+1}}(t-t_{x,i})(s-t_{x,j})\Phi^{(1,1)}(\epsilon_{x,i}^{(1)},\eta_{x,j}^{(1)})f(t)f(s)dtds.
\end{align*}
We first consider the term $I_{i,j}^{(1)}$. For $l=0,1,2$, let, \begin{equation}
\omega_{i,l}=\int_{t_{x,i}}^{t_{x,i+1}}(t-t_{x,i})^lf(t)dt \label{wildefin}
\end{equation} The term $I_{i,j}^{(1)}$ can then be written as,
\begin{equation}
I_{i,j}^{(1)}  =  \Phi^{(1,0)}(t_{x,i},t_{x,j}) \Big(  2d_{x,i} \omega_{i,0}\omega_{j,0}-4\omega_{i,1}\omega_{j,0} \Big).
\end{equation}
Expanding $f$ around $t_{x,i}$ yields, 
\begin{align}
\omega_{i,l} & = \int_{t_{x,i}}^{t_{x,i+1}}(t-t_{x,i})^l \big(f(t_{x,i}) + (t-t_{x,i}) f'(t_{x,i})+\frac{1}{2}(t-t_{x,i})^2 f''(\epsilon_{x,i}^{(4)}) \big) \;dt \nonumber\\
& = \frac{d_{x,i}^{(l+1)}}{(l+1)}f(t_{x,i})  + \frac{d_{x,i}^{(l+2)}}{(l+2)}f'(t_{x,i}) +O\Big(\frac{1}{n^{(l+3)}}\Big), \label{wil}
\end{align}
for some  $\epsilon_{x,i}^{(4)}$ in $]t_{x,i},t_{x,i+1}[$. Thus for $l=0,1,2$,
\begin{align*}
& I_{i,j}^{(1)}  = \Phi^{(1,0)}(t_{x,i},t_{x,j})  \Bigg( 2d_{x,i} \bigg( d_{x,i}f(t_{x,i})  + \frac{d_{x,i}^2}{2}f'(t_{x,i}) +O\Big(\frac{1}{n^3}\Big) \bigg) \\
&~~~ \times \bigg( d_{x,j}f(t_{x,j})  + \frac{d_{x,j}^2}{2}f'(t_{x,i}) +O\Big(\frac{1}{n^3}\Big)\bigg)\\
&-4\bigg(\frac{d_{x,i}^2}{2}f(t_{x,i})  + \frac{d_{x,i}^3}{3}f'(t_{x,i}) +O\Big(\frac{1}{n^4}\Big) \bigg)\bigg( d_{x,j}f(t_{x,j})  + \frac{d_{x,j}^2}{2}f'(t_{x,i}) +O\Big(\frac{1}{n^3}\Big)\bigg) \Bigg)\\
& = \Phi^{(1,0)}(t_{x,i},t_{x,j})\Big(-\frac{1}{3}f'(t_{x,i}) f(t_{x,j})d_{x,i}^3d_{x,j}+O\Big(\frac{1}{n^5}\Big)\Big).
\end{align*}
We obtain using Equations \eqref{O20O11} and  \eqref{d=1/nf},
\begin{align*}
I_{i,j}^{(1)} & = -\frac{1}{3}\Phi^{(1,0)}(t_{x,i},t_{x,j})f'(t_{x,i}) f(t_{x,j})d_{x,i}^3d_{x,j}+O\Big(\frac{1}{n^5h^3}\Big)\\
& = -\frac{1}{3n^2}\Phi^{(1,0)}(t_{x,i},t_{x,j})\frac{f'(t_{x,i})}{f^2(t_{x,i}^*)} f(t_{x,j})d_{x,i}d_{x,j}+O\Big(\frac{1}{n^5h^3}\Big),
\end{align*}
for some $t_{x,i}^*$ in $]t_{x,i},t_{x,i+1}[$. Using Lemma \ref{N_T=Onh} and the integrability of $\varphi_{x,h},\varphi_{x,h}',f,$ and of $R^{(0,1)}(.,t)$ and applying Lemma \ref{sumtointegral} twice,  we obtain,
\begin{align}
&\overset{{N_{T_n}-1}}{\underset{i \ne j=1}{\sum\sum}  }  I_{i,j}^{(1)} = -\frac{1}{3n^2} \overset{{N_{T_n}-1}}{\underset{i \ne j=1}{\sum\sum}  }  \Phi^{(1,0)}(t_{x,i},t_{x,j})\frac{f'(t_{x,i})}{f^2(t_{x,i}^*)} f(t_{x,j})d_{x,i}d_{x,j}+O\Big(\frac{1}{n^3h}\Big)\nonumber \\
& = -\frac{1}{3n^2}\int_{x-h}^{x+h}\int_{x-h}^{x+h} \Phi^{(1,0)}(t,s) \frac{ f'(t)}{f^2(t)}f(s)1_{\{s \ne t \}}\;dt\;ds +O\Big(\frac{1}{n^3h^2}\Big).\label{Iij1}
\end{align}
Similarly we verify that, 
\begin{align}
&\overset{{N_{T_n}-1}}{\underset{i \ne j=1}{\sum\sum}  }  I_{i,j}^{(2)} =-\frac{1}{3n^2}\int_{x-h}^{x+h}\int_{x-h}^{x+h} \Phi^{(0,1)}(t,s) \frac{ f'(s)}{f^2(s)}f(t)1_{\{s \ne t \}}\;dt\;ds +O\Big(\frac{1}{n^3h^2}\Big)\nonumber\\
& =-\frac{1}{3n^2}\int_{x-h}^{x+h}\int_{x-h}^{x+h} \Phi^{(1,0)}(t,s) \frac{ f'(t)}{f^2(t)}f(s)1_{\{s \ne t \}}\;dt\;ds +O\Big(\frac{1}{n^3h^2}\Big).\label{Iij2}
\end{align}
We now control the term $I_{i,j}^3$. We have,
\begin{align*}
&I_{i,j}^{(3)} =d_{x,i}^2 \Phi^{(2,0)}(t_{x,i},t_{x,j}) \int_{t_{x,i}}^{t_{x,i+1}}\int_{t_{x,j}}^{t_{x,j+1}}f(t)f(s)dtds\\
&-2\Phi^{(2,0)}(t_{x,i},t_{x,j})\int_{t_{x,i}}^{t_{x,i+1}}\int_{t_{x,j}}^{t_{x,j+1}}(t-t_{x,i})^2f(t)f(s)dtds\\
&+\frac{1}{2}d_{x,i}^2 \int_{t_{x,i}}^{t_{x,i+1}}\int_{t_{x,j}}^{t_{x,j+1}}\Phi^{(2,0)}(\epsilon_{x,i}^{(2)},t_{x,j})+ \Phi^{(2,0)}(\epsilon_{x,i}^{(3)},t_{x,j})-2\Phi^{(2,0)}(t_{x,i},t_{x,j})f(t)f(s)dtds\\
&-2\int_{t_{x,i}}^{t_{x,i+1}}\int_{t_{x,j}}^{t_{x,j+1}}(t-t_{x,i})^2\big(\Phi^{(2,0)}(\epsilon_{x,i}^{(1)},t_{x,j})-\Phi^{(2,0)}(t_{x,i},t_{x,j})\big)f(t)f(s)dtds.
\end{align*}
Using \eqref{O20O11}, Lemma \ref{N_T=Onh} and Equation \eqref{wildefin} we get,
\begin{align*}
I_{i,j}^{(3)}& =d_{x,i}^2 \Phi^{(2,0)}(t_{x,i},t_{x,j}) \omega_{i,0}\omega_{j,0}-2\Phi^{(2,0)}(t_{x,i},t_{x,j})\omega_{i,2}\omega_{i,0}+O\Big(\frac{1}{n^5h^5}\Big).
\end{align*}
 Note first that, using \eqref{wil} for $l=0$ along with $l=2$ and Lemma \ref{N_T=Onh}, we obtain,
\begin{align*}
I_{i,j}^{(3)}& = \frac{1}{3}\Phi^{(2,0)}(t_{x,i},t_{x,j})d_{x,i}^3d_{x,j}f(t_{x,i})f(t_{x,j})+O\Big(\frac{1}{n^5h^5}\Big)\\
& = \frac{1}{3n^2}\Phi^{(2,0)}(t_{x,i},t_{x,j})\frac{f(t_{x,i})}{f^2(t_{x,i}^*)}f(t_{x,j})d_{x,i}d_{x,j}+O\Big(\frac{1}{n^5h^5}\Big),
\end{align*}
Likewise, using Lemma \ref{N_T=Onh} and the integrability of $\varphi_{x,h}^{(k)},f^{(k)}$ for $k=0,1,2$ we have, 
\begin{align}
&\overset{{N_{T_n}-1}}{\underset{i \ne j=1}{\sum\sum}} I_{i,j}^{(3)} &= \frac{1}{3n^2}\int_{x-h}^{x+h}\int_{x-h}^{x+h} \Phi^{(2,0)}(t,s) \frac{f(s)}{f(t)}1_{\{s \ne t \}}\;dt\;ds +O\Big(\frac{1}{n^3h^3}\Big). \label{Iij3}
\end{align}
Similarly, we obtain,
\begin{align}
&\overset{{N_{T_n}-1}}{\underset{i \ne j=1}{\sum\sum}} I_{i,j}^{(4)} = \frac{1}{3n^2}\int_{x-h}^{x+h}\int_{x-h}^{x+h} \Phi^{(0,2)}(t,s) \frac{f(t)}{f(s)}1_{\{s \ne t \}}\;dt\;ds +O\Big(\frac{1}{n^3h^3}\Big)\nonumber \\
&= \frac{1}{3n^2}\int_{x-h}^{x+h}\int_{x-h}^{x+h} \Phi^{(2,0)}(t,s) \frac{f(s)}{f(t)}1_{\{s \ne t \}}\;dt\;ds +O\Big(\frac{1}{n^3h^3}\Big).\label{Iij4}
\end{align}
Finally, for the term $I_{i,j}^{(5)}$, we have,
\begin{align*}
& I_{i,j}^{(5)} = d_{x,i} d_{x,j}\Phi^{(1,1)}(t_{x,i},t_{x,j})\int_{t_{x,i}}^{t_{x,i+1}}\int_{t_{x,j}}^{t_{x,j+1}}f(t)f(s)dtds\\
&-4\Phi^{(1,1)}(t_{x,i},t_{x,j})\int_{t_{x,i}}^{t_{x,i+1}}\int_{t_{x,j}}^{t_{x,j+1}}(t-t_{x,i})(s-t_{x,j})f(t)f(s)dtds\\
&+ d_{x,i} d_{x,j}\int_{t_{x,i}}^{t_{x,i+1}}\int_{t_{x,j}}^{t_{x,j+1}}\Big(\Phi^{(1,1)}(\epsilon_{x,i}^{(3)},\eta_{x,j}^{(3)})-\Phi^{(1,1)}(t_{x,i},t_{x,j})\Big)f(t)f(s)dtds\\
&-4\int_{t_{x,i}}^{t_{x,i+1}}\int_{t_{x,j}}^{t_{x,j+1}}(t-t_{x,i})(s-t_{x,j})\Big(\Phi^{(1,1)}(\epsilon_{x,i}^{(1)},\eta_{x,j}^{(1)})-\Phi^{(1,1)}(t_{x,i},t_{x,j})\Big)f(t)f(s)dtds\\
&= d_{x,i} d_{x,j}\Phi^{(1,1)}(t_{x,i},t_{x,j})\omega_{i,0}\omega_{j,0}-4\Phi^{(1,1)}(t_{x,i},t_{x,j})\omega_{i,1}\omega_{j,1}\\
&+ d_{x,i} d_{x,j}\int_{t_{x,i}}^{t_{x,i+1}}\int_{t_{x,j}}^{t_{x,j+1}}\Big(\Phi^{(1,1)}(\epsilon_{x,i}^{(3)},\eta_{x,j}^{(3)})-\Phi^{(1,1)}(t_{x,i},t_{x,j})\Big)f(t)f(s)dtds\\
&-4\int_{t_{x,i}}^{t_{x,i+1}}\int_{t_{x,j}}^{t_{x,j+1}}(t-t_{x,i})(s-t_{x,j})\Big(\Phi^{(1,1)}(\epsilon_{x,i}^{(1)},\eta_{x,j}^{(1)})-\Phi^{(1,1)}(t_{x,i},t_{x,j})\Big)f(t)f(s)dtds.
\end{align*}
Recall that $f, f', \frac{1}{f} $ are all bounded and using \eqref{O20O11}  and \eqref{wil} with $l=l'=1$ we obtain,
\begin{equation*}
I_{i,j}^{(5)}=O\Big(\frac{1}{n^5h^5}\Big).
\end{equation*}
Finally, since $N_{T_n}=O(nh)$ from Lemma \ref{N_T=Onh}, we obtain,
\begin{equation}
\overset{{N_{T_n}-1}}{\underset{i \ne j=1}{\sum\sum}} I_{i,j}^{(5)}=O\Big(\frac{1}{n^3h^3}\Big). \label{Iij5}
\end{equation}
Replacing \eqref{Iij1}, \eqref{Iij2}, \eqref{Iij3}, \eqref{Iij4} and \eqref{Iij5} in \eqref{Iij} we obtain,

\begin{align*}
~&\overset{{N_{T_n}-1}}{\underset{i \ne j=1}{\sum\sum}} I_{i,j}
=\frac{1}{6mn^2}\int_{x-h}^{x+h}\int_{x-h}^{x+h} \bigg( \frac{\Phi^{(2,0)}(t,s)f(t)-\Phi^{(1,0)}(t,s)f'(t)}{f^2(t)}\bigg)1_{\{s \ne t \}} f(s)dtds\\
&~~~~~~~~+O\Big(\frac{1}{mn^3h^3}\Big)\\
& = \frac{1}{6mn^2}\int_{x-h}^{x+h} \bigg( \int_{x-h}^{s} \frac{\partial}{\partial s} \Big( \frac{\Phi^{(1,0)}(t,s)}{f(t)}\Big) \bigg)dt~f(s) ds +O\Big(\frac{1}{mn^3h^3}\Big)\\
&~~~~~~+\frac{1}{6mn^2}\int_{x-h}^{x+h} \bigg( \int_{s}^{x+h} \frac{\partial}{\partial s} \Big( \frac{\Phi^{(1,0)}(t,s)}{f(t)}\Big) \bigg)dt~f(s) ds +O\Big(\frac{1}{mn^3h^3}\Big)\\
& = \frac{1}{6mn^2}\int_{x-h}^{x+h} \big( \Phi^{(1,0)}(s^-,s)-\Phi^{(1,0)}(s^+,s)\big) ds +O\Big(\frac{1}{mn^3h^3}\Big)\\
&~~~~~~+\frac{1}{6mn^2}\int_{x-h}^{x+h} \bigg(\frac{\Phi^{(1,0)}(x+h,s)}{f(x+h)}-\frac{\Phi^{(1,0)}(x-h,s)}{f(x-h)} \bigg) f(s)~ds+O\Big(\frac{1}{mn^3h^3}\Big). 
\end{align*}
Note that for $t \ne s$, \begin{align}
\Phi^{(1,0)}(t,s) & = \bigg(\frac{\varphi_{x,h}'(t) f(t) -\varphi_{x,h}(t) f'(t)}{f^2(t)} R(t,s)+\frac{\varphi_{x,h}(t)}{f(t)}R^{(1,0)}(t,s) \bigg)\frac{\varphi_{x,h}(s)}{f(s)}.\label{expressionofphi}
\end{align}
It follows from \eqref{varphiis0} that,\[\frac{\Phi^{(1,0)}(x+h,s)}{f(x+h)}=\frac{\Phi^{(1,0)}(x-h,s)}{f(x-h)} = 0 ~~~~\text{for all}~ s\in]x-h,x+h[.\]
Thus,
\begin{equation}
\overset{{N_{T_n}-1}}{\underset{i \ne j=1}{\sum\sum}} I_{i,j}=\frac{1}{6mn^2}\int_{x-h}^{x+h} \big( \Phi^{(1,0)}(t^-,t)-\Phi^{(1,0)}(t^+,t)\big) dt+O\Big(\frac{1}{mn^3h^3}\Big). \label{offdiagtrap}
\end{equation}
Inserting \eqref{diagtermtrap} and \eqref{offdiagtrap} in \eqref{var-sigma2}, we obtain,
{\small
\begin{align}
\Var (\hat{g}^{trap}_n (x) ) &= \frac{1}{m}\sigma_{x,h}^2 + \frac{1}{12mn^2} \int_{x-h}^{x+h} \big( \Phi^{(1,0)}(t^-,t)-\Phi^{(1,0)}(t^+,t)\big)\;dt\nonumber\\
&~~~~~~~~+ O\Big(\frac{1}{mn^3h^3}\Big). \label{var=danspreuve}
\end{align}}
Applying \eqref{expressionofphi} it follows that,
\begin{align}
\Phi^{(1,0)}(t^-,t)-\Phi^{(1,0)}(t^+,t) & = \frac{\varphi_{x,h}^2(t)}{f^2(t)}\Big(R^{(1,0)}(t^-,t)-R^{(1,0)}(t^-,t)\Big) = \frac{\varphi_{x,h}^2(t)}{f^2(t)} \alpha(t).\label{secondone}
\end{align}
Replacing  \eqref{secondone} in \eqref{var=danspreuve} we obtain, 
\begin{align}
\Var (\hat{g}^{trap}_n (x) ) = \frac{1}{m}\sigma_{x,h}^2 + \frac{1}{12mn^2} \int_{x-h}^{x+h} \frac{\varphi_{x,h}^2(t)}{f^2(t)} \alpha(t)\;dt+ O\Big(\frac{1}{mn^3h^3}\Big). \label{mpol}
\end{align}
Since $\alpha$ and $f$ are continuous on $[0,1]$, then one can write,
\begin{align}
&\int_{x-h}^{x+h}\frac{\alpha(t)}{f^2(t)} \varphi_{x,h}^{2}(t)dt  = \frac{1}{h}\int_{-1}^1 \frac{\alpha(x-th)}{f^2(x-th)} K^2(t)\;dt\nonumber\\
& = \frac{1}{h}\frac{\alpha(x)}{f^2(x)}\int_{-1}^1  K^2(t)\;dt+\frac{1}{h}\int_{-1}^1 \bigg(\frac{\alpha(x-th)}{f^2(x-th)} -\frac{\alpha(x)}{f^2(x)}\bigg) K^2(t)\;dt\nonumber\\
& = \frac{1}{h}\frac{\alpha(x)}{f^2(x)}\int_{-1}^1  K^2(t)\;dt + O(1). \label{lmpl}
\end{align}
Recall that for an even kernel, we have a simplified expression of $\sigma_{x,h}^2$ given by Benhenni and Rachdi \cite{Benhenni and Rachdi} as follows,
\begin{equation}
\sigma_{x,h}^2 =  R(x,x)-\frac{1}{2}\alpha(x)C_K h  + o(h),  \label{expansionofsigma}
\end{equation}
where $C_K = \int_{-1}^1 \int_{-1}^1 |u-v|K(u)K(v) du dv.$ \\
Finally, using \eqref{lmpl} and \eqref{expansionofsigma} in \eqref{mpol} yields,
\begin{align*}
\Var (\hat{g}^{trap}_n (x) ) & =  \frac{1}{m}\Big( R(x,x)-\frac{1}{2}\alpha(x)C_K h \Big)+\frac{1}{12mn^2h} \frac{\alpha(x)}{f^2(x)}\int_{-1}^1  K^2(t)\;dt \\
&~~~~~~+ o\Big(\frac{h}{m}\Big)+ O\Big(\frac{1}{mn^2}+\frac{1}{mn^3h^3}\Big) . 
\end{align*}
This concludes the proof of Proposition \ref{variancelemma1}. $\Box$
\subsection{Proof of Proposition \ref{optimalh}.}
Let $I_1=\int_{0}^{1}R(x,x)w(x)\;dx$, $ I_2=\int_{0}^{1}  \frac{\alpha(x)}{f^2(x)} w(x)\;dx$ and put, \[\Psi (h,m)=-\frac{C_K h}{2m} \int_{0}^{1}\alpha(x)  w(x)\;dx + \frac{1}{4} h^4 B^2 \int_{0}^{1}[g''(x)]^2 w(x)\;dx.\]
We have from Equation \eqref{IMSEexpression} in Theorem \ref{IMSE},
{\small
\begin{align*}
\IMSE (h) & = \frac{I_1}{m} +\Psi(h,m) + \frac{VI_2}{12mn^2h} +o \Big(h^4+\frac{h}{m}\Big)+O\Big(\frac{1}{n^3h}+\frac{1}{mn^3h^3}+\frac{1}{mn^2}+\frac{1}{n^6h^6}\Big).
\end{align*}}
\noindent Let $h^*$ be as defined in \eqref{hoptimal}. It is clear that  $h^* = \underset{0< h < 1}{\argmin} ~\Psi (h,m)$ so that $\Psi (h,m) \geq \Psi (h^*,m)$  for every $0< h < 1$. Let $h_{n,m}$ be as defined in Corollary \ref{optimalh}. We have,
{\small
\begin{align*}
&~~~~~~~~~~~~~~~~~~~~~~~~~~~~~~~~~~~~~~~~~~\frac{\IMSE(h^*)}{\IMSE(h_{n,m})}\\
 &= \frac{ \frac{I_1}{m}  +\Psi(h^*,m) + \frac{VI_2}{12mn^2h^*} +o \Big({h^*}^4+\frac{h^*}{m}\Big)+O\Big(\frac{1}{n^3{h^*}}+\frac{1}{mn^3{h^*}^3}+\frac{1}{mn^2}+\frac{1}{n^6{h^*}^6}\Big)}{\frac{I_1}{m}  +\Psi(h_{n,m},m) + \frac{VI_2}{12mn^2h_{n,m}} +o \Big(h_{n,m}^4+\frac{h_{n,m}}{m}\Big)+O\Big(\frac{1}{n^3h_{n,m}}+\frac{1}{mn^3h_{n,m}^3}+\frac{1}{mn^2}+\frac{1}{n^6h_{n,m}^6}\Big)}\\
&\\
& \leq \frac{ I_1  +m\Psi(h_{n,m},m) + \frac{VI_2}{12n^2h^*}+o \Big({mh^*}^4+h^*\Big)+O\Big(\frac{m}{n^3{h^*}}+\frac{1}{n^3{h^*}^3}+\frac{1}{n^2}+\frac{m}{n^6{h^*}^6}\Big)}{I_1  +m\Psi(h_{n,m},m) + \frac{VI_2}{12n^2h_{n,m}}  +o \Big(mh_{n,m}^4+h_{n,m}\Big)+O\Big(\frac{m}{n^3h_{n,m}}+\frac{1}{n^3h_{n,m}^3}+\frac{1}{n^2}+\frac{m}{n^6h_{n,m}^6}\Big)}.
\end{align*}}
\noindent  Using the definition of $h^*,$ $mh_{n,m}^3=O(1)$,  $\underset{n,m \to \infty}{\lim}~h_{n,m}=0$  and the assumption $\frac{m}{n}=O(1)$ as $n,m \to \infty$ we know that $m\Psi(h_{n,m},m)=O(h_{n,m})$. Then,  \[ \underset{n,m \to \infty}{\overline {\lim}} \frac{\IMSE(h^*)}{\IMSE(h_{n,m})} \leq 1.\]This concludes the proof of Proposition \ref{optimalh}. $\Box$
\subsection{Proof of Corollary \ref{optimaldensitytheorem}.}
Let $f^*$ be as defined in \eqref{f^*}. Let $D(f) = \int_{0}^{1} \frac{\alpha(x)}{f^2(x)}w(x) \;dx,$
then it is sufficient to prove that: \[D(f^*)\leq D(f)~~~ \text{for every positive density}~ f~\text{ on } [0,1]. \] 
Applying Hölder's inequality, we get, 
\begin{align*}
D(f^*) & = \bigg(\int_{0}^{1} \{\alpha(x) w(x) \}^{1/3} \;dx\bigg)^3 = \bigg (\int_{0}^{1} \Big(\frac{\alpha(x) w(x)}{f^2(x)} \Big)^{1/3} f^{2/3}(x) \;dx\bigg)^3\\
& \leq \bigg (\int_{0}^{1}\frac{\alpha(x) w(x)}{f^2(x)}  \;dx\bigg)\bigg (\int_{0}^{1} f(x)\;dx \bigg)^2 = D(f).
\end{align*}
Hence, \[\underset{\{f>0 ~\text{density on}~ [0,1] \}}{\text{argmin}} D(f)=f^*. \]
This completes the proof of Corollary \ref{optimaldensitytheorem}. $\Box$
\subsection{Proof of Theorem \ref{minimaxtheorem}.}
Let $f^*$ be as defined in \eqref{f^*}. The proof of this theorem will be done in two steps:
\begin{enumerate}
\item $\sup\{\Psi_{(\alpha,w)}(f^*) / (\alpha, w) \in \Lambda \} \leq \epsilon_1\epsilon_2$.
\item $\forall f, \exists (\alpha,w) \in \Lambda : \Psi_{(\alpha,w)}(f) \geq \epsilon_1\epsilon_2$.
\end{enumerate}
\textbf{First step:}  By direct application of the Hölder's inequality we have: 
\begin{align*}
\Psi_{(\alpha,w)}(f^*) & = \bigg(\int_{0}^{1} \{\alpha(s) w(s) \}^{1/3} \;ds\bigg)^3 =  \bigg(\int_{0}^{1} \alpha(s)^{1/3} \sqrt{w(s)}^{2/3}  \;ds\bigg)^3 \\
& \leq \Big (  \int_{0}^{1} \alpha(s)ds \Big) \Big ( \int_{0}^{1} \sqrt{w(s)}ds \Big)^2 \leq \epsilon_1 \epsilon_2.
\end{align*}
\textbf{Second step:}  Let $f$ be an arbitrary positive density. Take $\alpha^*\equiv\epsilon_1$ and $w^*\equiv \epsilon_2$, then $(\alpha^*,w^*) \in \Lambda$ and: 
$$\Psi_{(\alpha^*,w^*)}(f) = \int_{0}^{1} \frac{\alpha^*(s)w^*(s)}{f^2(s)}~ds=\epsilon_1 \epsilon_2\int_{0}^{1} \frac{1}{f^2(s)}~ds \geq \epsilon_1 \epsilon_2\ ,$$
since, using the Hölder's inequality we have:
\begin{align*}
1&= \int_{0}^{1} f^{2/3}(s)\Big(\frac{1}{f^{2}(s)}\Big)^{1/3}~ds \leq \Big(\int_{0}^{1} f(s)~ds\Big)^{2/3} \Big(\int_{0}^{1}\frac{1}{f^2(s)}~ds\Big)^{1/3}= \Big(\int_{0}^{1}\frac{1}{f^2(s)}~ds\Big)^{1/3}.
\end{align*}
This completes the proof of Theorem \ref{minimaxtheorem}. $\Box$
\subsection{Proof of Theorem \ref{asymptotic normality theorem}.}
Let $x \in ]0,1[$ be fixed. We have, 
{\small
\begin{equation} \label{diffasym}
\sqrt{m} \Big(\hat{g}_{n,m}^{trap}(x) -g(x) \Big)=\sqrt{m} \Big(\hat{g}_{n,m}^{trap}(x) -\mathbb{E}\big( \hat{g}_{n,m}^{trap}(x)\big) \Big)+\sqrt{m}~\Bias \Big(\hat{g}_{n,m}^{trap}(x)\Big).
\end{equation}}
Since $\underset{n,m \to \infty}{\lim} \sqrt{m} h^2 =0$ and $\underset{n,m \to \infty}{\lim} nh^2 =\infty$ then Proposition \ref{biaslemma} implies that,
\begin{equation}
\underset{n,m \to \infty}{\lim}\sqrt{m}~\Bias \Big(\hat{g}_{n,m}^{trap}(x)\Big) =0.
\end{equation}
Consider now the first term of the right side of \eqref{diffasym}. Since $\overline{Y}(t_{x,i})-\mathbb{E}(\overline{Y}(t_{x,i})) =\overline{\varepsilon}(t_{x,i})$, we have, as done by Fraiman and Pérez Iribarren  \cite{Fraiman},
{\small
\begin{align}
& \sqrt{m} \Big(\hat{g}_{n,m}^{trap}(x) -\mathbb{E}\big( \hat{g}_{n,m}^{trap}(x)\big) \Big)  =  \frac{1}{\sqrt{m}} \Big\{  \sum_{j=1}^{m} \frac{1}{2n}  \sum_{i=1}^{N_{T_n}-1}  \Big( \big(\frac{\varphi_{x,h}}{f} \varepsilon_j \big) (t_{x,i} ) +\big(\frac{\varphi_{x,h}}{f} \varepsilon_j \big) (t_{x,i+1} )  \Big) \Big\} \nonumber\\
& =  \frac{1}{\sqrt{m}} \sum_{j=1}^{m} \frac{1}{2n}\sum_{i=1}^{N_{T_n}-1}\frac{\varphi_{x,h}}{f}(t_{x,i} )  \big (\varepsilon_j (t_{x,i} )-\varepsilon_j(x) \big)\nonumber\\
& ~~~~~+ \frac{1}{\sqrt{m}} \sum_{j=1}^{m} \frac{1}{2n}\sum_{i=1}^{N_{T_n}-1}\frac{\varphi_{x,h}}{f}(t_{x,i+1} )  \big (\varepsilon_j (t_{x,i+1} )-\varepsilon_j(x) \big) \nonumber\\
& ~~~~~+\bigg( \frac{1}{2n}  \sum_{i=1}^{N_{T_n}-1} \Big( \frac{\varphi_{x,h}}{f}(t_{x,i} ) +\frac{\varphi_{x,h}}{f} (t_{x,i+1} ) \Big) \bigg) \bigg(\frac{1}{\sqrt{m}} \sum_{j=1}^{m} \varepsilon_j(x)\bigg). \label{sw}
\end{align}}
\noindent We start by controlling the last term of this last equation.  Recall that Equation \eqref{d=1/nf} yields for some $t_{x,i}^* \in ]t_{x,i},t_{x,i+1}[$ that $\frac{1}{n}=(t_{x,i+1}-t_{x,i})f(t_{x,i}^*)$. From the Riemann integrability of $\varphi_{x,h}$ and $f$ and Lemma \ref{sumtointegral}  we obtain,
\begin{align*}
&~~~~~~~~\frac{1}{2n}  \sum_{i=1}^{N_{T_n}-1} \Big( \frac{\varphi_{x,h}}{f}(t_{x,i} ) +\frac{\varphi_{x,h}}{f} (t_{x,i+1} ) \Big)= \\
&  \frac{1}{2}  \sum_{i=1}^{N_{T_n}-1} \Big( \frac{\varphi_{x,h}}{f}(t_{x,i} ) +\frac{\varphi_{x,h}}{f} (t_{x,i+1} ) \Big)f(t_{x,i}^*)(t_{x,i+1}-t_{x,i}) \underset{m,n \to \infty}{\longrightarrow} \int_{-1}^{1}K(t)~dt= 1.
\end{align*}
where $d_{x,i}=t_{x,i+1}-t_{x,i}$ and $t_{x,i}^* \in ]t_{x,i},t_{x,i+1}[$. 
The Central Limit Theorem for i.i.d. variables yields,
\[ \frac{1}{\sqrt{m}} \sum_{j=1}^{m} \varepsilon_j(x) \underset{m \to \infty}{\overset{\mathscr D}{\longrightarrow}} Z ~~~where~~ Z \sim \mathcal{N}(0,R(x,x)). \]
We shall prove now that the two first terms of Equation \eqref{sw} tend to 0 in probability as $n,m$ tends to infinity. We will only study the first term, the second one is treated analogously. Let,
\begin{align*}
A_{m,n}(x) & = \frac{1}{\sqrt{m}} \sum_{j=1}^{m} \frac{1}{2n}\sum_{i=1}^{N_{T_n}-1}\frac{\varphi_{x,h}}{f}(t_{x,i} )  \big (\varepsilon_j (t_{x,i} )-\varepsilon_j(x) \big)\overset{\Delta}{=}\frac{1}{\sqrt{m}} \sum_{j=1}^{m} T_{n,j}(x).
\end{align*}
From the Chebyshev inequality, it suffices to prove that $\underset{n,m \to \infty}{\lim} \mathbb{E}(A_{m,n}^2(x)) =0 $. We have for $j \ne l$,  $\mathbb{E} (\varepsilon_j(x)\varepsilon_l(y))=0$ so $\mathbb{E} (T_{n,j}(x) T_{n,l}(x)) =0$. Hence,
\begin{align*}
\mathbb{E}(A_{m,n}^2(x)) & = \frac{1}{m} \sum_{j=1}^{m} \sum_{l=1}^{m} \mathbb{E} (T_{n,j}(x) T_{n,l}(x)) = \frac{1}{m} \sum_{j=1}^{m} \mathbb{E} (T_{n,j}^2(x)).
\end{align*} 
We have, 
\begin{align*}
&~~~~~~~~~~~~~~\mathbb{E} (T_{n,j}^2(x))  =\\
&\frac{1}{4n^2} \sum_{i=1}^{N_{T_n}-1}\sum_{k=1}^{N_{T_n}-1}\frac{\varphi_{x,h}}{f}(t_{x,i} )\frac{\varphi_{x,h}}{f}(t_{x,k} ) \mathbb{E} \Big(  \big (\varepsilon_j (t_{x,i} )-\varepsilon_j(x) \big) \big (\varepsilon_j (t_{x,k} )-\varepsilon_j(x) \big) \Big)\\
& = \frac{1}{4n^2} \sum_{i=1}^{N_{T_n}-1}\sum_{k=1}^{N_{T_n}-1}\frac{\varphi_{x,h}}{f}(t_{x,i} )\frac{\varphi_{x,h}}{f}(t_{x,k} ) \Big(R(t_{x,i},t_{x,k})-R(t_{x,i},x)-R(x,t_{x,k})+R(x,x)\Big).
\end{align*} 
Since $\mathbb{E} ((T_{n,j}^2(x))$ does not depend on $j$ we get, 
\begin{align}
&~~~~~~~~~~~~~~\mathbb{E}(A_{m,n}^2(x))= \nonumber\\
& \frac{1}{4n^2} \sum_{i=1}^{N_{T_n}-1}\sum_{k=1}^{N_{T_n}-1}\frac{\varphi_{x,h}}{f}(t_{x,i} )\frac{\varphi_{x,h}}{f}(t_{x,k} ) \Big(R(t_{x,i},t_{x,k})-R(t_{x,i},x)-R(x,t_{x,k})+R(x,x)\Big)\nonumber\\
& \overset{\Delta}{=} \frac{1}{4} \Big(B_{n,1}(x)-B_{n,2}(x)-B_{n,3}(x)+B_{n,4}(x)\Big).\label{sumBn}
\end{align}
We obtain using Equation \eqref{d=1/nf} for $t_{x,i}^* \in ]t_{x,i},t_{x,i+1}[$,
\begin{align*}
B_{n,1}(x) & =  \sum_{i=1}^{N_{T_n}-1}\sum_{k=1}^{N_{T_n}-1} f(t_{x,i}^*)f(t_{x,k}^*)~\frac{\varphi_{x,h}}{f}(t_{x,i} )\frac{\varphi_{x,h}}{f}(t_{x,k} ) R(t_{x,i},t_{x,k})d_{x,i}d_{x,k}.
\end{align*}
The use of Lemma \ref{sumtointegral} twice yields,
\begin{align*}
B_{n,1}(x) & =  \sum_{i=1}^{N_{T_n}-1} f(t_{x,i}^*)\frac{\varphi_{x,h}}{f}(t_{x,i} )d_{x,i} \Big \{ \int_{x-h}^{x+h}\varphi_{x,h}(t ) R(t_{x,i},t)~ dt +O(\frac{1}{nh}) \Big\} \\
& =  \int_{x-h}^{x+h}\varphi_{x,h}(t )   \Big \{  \sum_{i=1}^{N_{T_n}-1} f(t_{x,i}^*)\frac{\varphi_{x,h}}{f}(t_{x,i} ) R(t_{x,i},t)  d_{x,i}  \Big\} dt+O(\frac{1}{nh})\\
& = \int_{x-h}^{x+h} \int_{x-h}^{x+h} \varphi_{x,h}(s)\varphi_{x,h}(t)R(s,t)~ds~dt+ O(\frac{1}{nh})=\sigma_{x,h}^2+ O(\frac{1}{nh}).
\end{align*}
Using \eqref{expansionofsigma} we obtain, \[B_{n,1}(x)  = R(x,x)-\frac{1}{2}\alpha(x)C_K h  + o(h)+ O(\frac{1}{nh}).\]
where $C_K = \int_{-1}^1 \int_{-1}^1 |u-v|K(u)K(v) du dv.$ Since $\underset{n \to \infty}{\lim}h=0$ and $\underset{n \to \infty}{\lim}nh=\infty$. Thus, 
\begin{equation}\label{Bn1}
\underset{n \to \infty}{\lim} B_{n,1}(x) = R(x,x).
\end{equation}
Consider now the term $B_{n,2}(x)$. We obtain using Lemma \ref{sumtointegral} twice,
\begin{align*}
B_{n,2}(x) & = \int_{x-h}^{x+h} \int_{x-h}^{x+h} \varphi_{x,h}(s)\varphi_{x,h}(t)R(s,x)~ds~dt+ O(\frac{1}{nh})\\
&= \int_{x-h}^{x+h}  \varphi_{x,h}(s)R(s,x)~ds+ O(\frac{1}{nh})\\
& = \int_{-1}^{1}  K(s)R(x-hs,x)~ds+ O(\frac{1}{nh})\\
& = \int_{-1}^{0}  K(s)R(x-hs,x)~ds+\int_{0}^{1}  K(s)R(x-hs,x)~ds+  O(\frac{1}{nh}).
\end{align*}
For $s \in ]-1,0[$, Taylor expansion of $R(\cdot,x)$ around $x$ yields, \[R(s,x)  = R(x-sh,x)-shR^{(1,0)}(x+,x)+o(h).\]
Similarly for  $s \in ]0,1[$ we obtain,\[R(x-sh,x)  = R(x,x)-shR^{(1,0)}(x-,x)+o(h).\]
Thus,
\begin{align*}
B_{n,2}(x) & = R(x,x)-hR^{(1,0)}(x+,x)\int_{-1}^{0}s~K(s)~ds\\
&~~-hR^{(1,0)}(x-,x)\int_{0}^{1}s~K(s)~ds+o(h)+ O(\frac{1}{nh}).
\end{align*}
Hence, \begin{equation} \label{Bn2}
\underset{n \to \infty}{\lim} B_{n,2}(x) = R(x,x).
\end{equation}
Similarly, \begin{equation}\label{Bn3}
\underset{n \to \infty}{\lim} B_{n,3}(x) = R(x,x).
\end{equation}
It is easy to see that,\begin{align}\label{Bn4}
\underset{n \to \infty}{\lim} B_{n,4}(x)& = \underset{n \to \infty}{\lim}  R(x,x)\sum_{i=1}^{N_{T_n}-1}\sum_{k=1}^{N_{T_n}-1}\frac{\varphi_{x,h}}{f}(t_{x,i} )\frac{\varphi_{x,h}}{f}(t_{x,k}) \nonumber)\\
& = R(x,x) \Big(\int_{-1}^{1}K(t)~dt\Big)^2=R(x,x).
\end{align}
Inserting \eqref{Bn1}, \eqref{Bn2}, \eqref{Bn3} and \eqref{Bn4} in \eqref{sumBn} yields, \[\underset{n,m \to \infty}{\lim} \mathbb{E}(A_{m,n}^2(x)) =0. \]
This concludes the proof of Theorem \ref{asymptotic normality theorem}. $\Box$

\section*{Appendix}
\begin{lemma}[Integral approximation of a sum] \label{sumtointegral}
	Let $u$ and $v$ be two Lipschitz functions  on $[x-h,x+h]$, i.e, there exists two positive numbers $l_1$ and  $l_2$ such that, \[ |u(s)-u(t)| \leq l_1|s-t|,~~~~|v(s)-v(t)| \leq l_2|s-t|. \] 
	Let $t_{x,1}< \cdots < t_{x,N_{T_n}}$ be points in $[x-h,x+h]$ and put $d_{x,i}=t_{x,i+1}-t_{x,i}$. Then,
	\[\sum_{i=1}^{N_{T_n}-1} u(t_{x,i}) v(t_{x,i}') d_{x,i} = \int_{x-h}^{x+h} u(t) v(t)~ dt +\Delta_{n,h},\]
	for any $t_{x,i}' \in [t_{x,i},t_{x,i+1}]$ for all $i=1,\cdots,n$ and for some appropriate positive constants $c_1,c_2$ and $c_3$,
	\[ |\Delta_{n,h}|\leq  c_1~ l_1 \frac{h}{n}~ \underset{t \in [0,1]}{\sup}|v(t)|~+c_2~ l_2 \frac{h}{n}~ \underset{t \in [0,1]}{\sup}|u(t)|+2\frac{c_3}{n} \underset{\underset{\cup [t_{x,\tiny {N_{T_n}}},x+h]}{t \in [x-h,t_{x,1}] }}{\sup}|v(t)u(t)| .\]
\end{lemma}
\textbf{Proof of Lemma \ref{sumtointegral}.}
In fact, let $\Delta_{x,h}=A-B$ where, \[ A  = \sum_{i=1}^{N_{T_n}-1} u(t_{x,i})v(t_{x,i}') d_{x,i} ~~~\text{and}~~ B =  \int_{x-h}^{x+h} u(t) v(t)~ dt. \]
We have, 
\begin{align*}
B & =   \sum_{i=1}^{N_{T_n}-1}  \int_{t_{x,i}}^{t_{x,i+1}} u(t) v(t)~ dt +\int_{x-h}^{t_{x,1}} u(t) v(t)~ dt+\int_{t_{x,N_{T_n}}}^{x+h} u(t) v(t)~ dt  \overset{\Delta}{=} B_1+B_2,
\end{align*}
where $B_2=\int_{x-h}^{t_{x,1}} u(t) v(t)~ dt+\int_{t_{x,N_{T_n}}}^{x+h} u(t) v(t)~ dt $.
On the one hand, since $(t_{x,1}-(x-h)) \leq  \underset{1 \leq i \leq n}{\sup}~d_{x,i}$ and $ (x+h-t_{x,N_{T_n}}) \leq  \underset{1 \leq i \leq n}{\sup}~d_{x,i}$ we have,\[| B_2 | \leq 2 c_3 \underset{\underset{\cup [t_{x,N_{T_n}},x+h]}{t \in [x-h,t_{x,1}] }}{\sup}|v(t)u(t)|   \underset{1 \leq i \leq n}{\sup}~d_{x,i}.\]
On the other hand, we have,
\begin{align*}
A- B_1 & =  \sum_{i=1}^{N_{T_n}-1}  \int_{t_{x,i}}^{t_{x,i+1}}\big( u(t_{x,i}) v(t_{x,i}')-  u(t) v(t)\big)~ dt\\
& =  \sum_{i=1}^{N_{T_n}-1}   v(t_{x,i}') \int_{t_{x,i}}^{t_{x,i+1}} \big( u(t_{x,i})-  u(t) \big)~ dt + \sum_{i=1}^{N_{T_n}-1}  \int_{t_{x,i}}^{t_{x,i+1}} u(t)  \big( v(t_{x,i}')- v(t)\big)~ dt.
\end{align*}
Since $u$ and $v$ are Lipschitz continuous we obtain,
\begin{align*}
|A- B_1| & \leq N_{T_n} \underset{t \in [0,1]}{\sup} |v(t)| l_1  \underset{1 \leq i \leq n}{\sup}~d_{x,i}^2  + N_{T_n} \underset{t \in [0,1]}{\sup} |u(t)| l_2  \underset{1 \leq i \leq n}{\sup}~d_{x,i}^2.
\end{align*}
Since $nh \geq 1$, Lemma \ref{N_T=Onh} yields that $\underset{1 \leq i \leq n}{\sup}~d_{x,i} = O(\frac{1}{n})$ and $N_{T_n}=O(nh)$. Hence,
\begin{align*}
|\Delta_{n,h}| & = |A- B| \leq |A-B_1|+|B_2|\\
&  \leq  c_1~ l_1 \frac{h}{n}~ \underset{t \in [0,1]}{\sup}|v(t)|~+c_2~ l_2 \frac{h}{n}~ \underset{t \in [0,1]}{\sup}|u(t)|+2\frac{c_3}{n} \underset{\underset{\cup [t_{x,\tiny {N_{T_n}}},x+h]}{t \in [x-h,t_{x,1}] }}{\sup}|v(t)u(t)|.
\end{align*}
This concludes the proof of Lemma \ref{sumtointegral}. $\Box$
\end{document}